\documentclass[12pt]{article}

\usepackage{amsfonts,amsmath,amssymb,amsthm}
\usepackage{bbm,booktabs,eurosym,fancyhdr,graphicx,multirow,setspace,url,xcolor}
\usepackage[T1]{fontenc}

\begin{document}

\newcommand{\real}{{\ensuremath{\mathbb{R}}}}
\newcommand{\dd}{{\ensuremath{\rm d}}}

\newcommand{\cA}{{\ensuremath{\cal A}}}
\newcommand{\cC}{{\ensuremath{\cal C}}}
\newcommand{\cF}{{\ensuremath{\cal F}}}
\newcommand{\cI}{{\ensuremath{\cal I}}}
\newcommand{\cN}{{\ensuremath{\cal N}}}
\newcommand{\cP}{{\ensuremath{\cal P}}}
\newcommand{\cS}{{\ensuremath{\cal S}}}

\newcommand{\myS}{\mbox{\rm S}}
\newcommand{\mySB}{\mbox{\rm S}_\theta^{\rm B}}
\newcommand{\mySE}{\mbox{\rm S}_{\alpha,\theta}^{\rm E}}
\newcommand{\mySQ}{\mbox{\rm S}_{\alpha,\theta}^{\rm Q}}

\newcommand{\myCE}{\cS_{\alpha}^{\rm E}}
\newcommand{\myCQ}{\cS_{\alpha}^{\rm Q}}

\newcommand{\myT}{{\ensuremath{\rm T}}}

\newcommand{\hsp}{\hspace{0.2mm}}
\newcommand{\vsp}{\vspace{-2mm}}

\newcommand{\myE}{{\ensuremath{\mathbb{E}}}}
\newcommand{\myQ}{{\ensuremath{\mathbb{Q}}}}
\newcommand{\one}{\ensuremath{\mathbbm{1}}}

\newcommand{\bysame}{\raisebox{1.5mm}{\underline{\hspace{2.75em}}}}
\newcommand{\done}{\hfill $\Box$}
\newcommand{\marginal}[1]{\marginpar{\raggedright\scriptsize #1}}

\newcommand{\beq}{\begin{equation}}
\newcommand{\eeq}{\end{equation}}
\newcommand{\beqa}{\begin{eqnarray}}
\newcommand{\eeqa}{\end{eqnarray}}
\newcommand{\beqas}{\begin{eqnarray*}}
\newcommand{\eeqas}{\end{eqnarray*}}

\newcommand{\supnl}{{\sup\nolimits}}

\newcommand{\ceq}{\!\!\! & = & \!\!\!}

\newcommand{\tcr}{\color{red}}
\newcommand{\tcb}{\color{blue}}

\pagestyle{plain}


\begin{center}
{

\LARGE \bf 

Of Quantiles and Expectiles: \\

Consistent Scoring Functions, Choquet Representations, and 
Forecast Rankings 

} 

\bigskip
Werner Ehm$^1$, Tilmann Gneiting$^{1,2}$, Alexander Jordan$^{1,2}$, Fabian Kr\"uger$^1$ 

\bigskip
{
\footnotesize 
$^1$Heidelberger Institut f\"ur Theoretische Studien, Computational Statistics Group \\
$^2$Karlsruher Institut f\"ur Technologie, Institut f\"ur Stochastik
}

\bigskip
\today
\end{center}

\vspace{2mm}
\begin{abstract}

In the practice of point prediction, it is desirable that forecasters
receive a directive in the form of a statistical functional, such as
the mean or a quantile of the predictive distribution.  When
evaluating and comparing competing forecasts, it is then critical that
the scoring function used for these purposes be consistent for the
functional at hand, in the sense that the expected score is minimized
when following the directive.

We show that any scoring function that is consistent for a quantile or
an expectile functional, respectively, can be represented as a mixture
of extremal scoring functions that form a linearly parameterized
family.  Scoring functions for the mean value and probability
forecasts of binary events constitute important examples.  The
quantile and expectile functionals along with the respective extremal
scoring functions admit appealing economic interpretations in terms of
thresholds in decision making.

The Choquet type mixture representations give rise to simple checks of
whether a forecast dominates another in the sense that it is
preferable under any consistent scoring function.  In empirical
settings it suffices to compare the average scores for only a finite
number of extremal elements.  Plots of the average scores with respect
to the extremal scoring functions, which we call Murphy diagrams,
permit detailed comparisons of the relative merits of competing
forecasts.

\bigskip
\noindent


\end{abstract}

\newpage

\smallskip
\section{Introduction}  \label{sec:introduction}

Over the past two decades, a broad transdisciplinary consensus has
developed that forecasts ought to be probabilistic in nature, i.e.,
they ought to take the form of predictive probability distributions
over future quantities or events (Gneiting and Katzfuss 2014).
Nevertheless, a wealth of applied settings require point forecasts, be
it for reasons of decision making, tradition, reporting requirements,
or ease of communication.  In this situation, a directive is required
as to the specific feature or functional of the predictive
distribution that is being sought.

We follow Gneiting (2011) and consider a functional to be a
potentially set-valued mapping $\myT(F)$ from a class of probability
distributions, $\cF$, to the real line, $\real$, with the mean or
expectation functional, quantiles, and expectiles being key examples.
Competing point forecasts are then compared by using a nonnegative
scoring function, $\myS(x,y)$, that represents the loss or penalty
when the point forecast $x$ is issued and the observation $y$
realizes.  A critically important requirement on the scoring function
is that it be consistent for the functional $\myT$ relative to the
class $\cF$, in the sense that
\begin{equation}  \label{eq:consistent.intro} 
\myE_F \hsp [\myS(t,Y)] \leq \myE_F \hsp [\myS(x,Y)] 
\end{equation} 
for all probability distributions $F \in \cF$, all $t \in \myT(F)$,
and all $x \in \real$.  If equality in (\ref{eq:consistent.intro}) 
implies that $x \in \myT(F)$, then the scoring function is strictly consistent.

To give a prominent example, the ubiquitous squared error scoring
function, $\myS(x,y) = (x-y)^2$, is strictly consistent for the mean
or expectation functional relative to the class of probability
distributions with finite variance.  However, there are many
alternatives.  In a classical paper, Savage (1971) showed, subject to
weak regularity conditions, that a scoring function is consistent for
the mean functional if and only if it is of the form
\begin{equation}  \label{eq:Bregman} 
\myS(x,y) = \phi(y) - \phi(x) -\phi'(x) \hsp (y-x), 
\end{equation} 
where the function $\phi$ is convex with subgradient $\phi'$; squared
error arises when $\phi(t) = t^2$.  Holzmann and Eulert (2014) proved
that when forecasts make ideal use of nested information bases, the
forecast with the broader information basis is preferable under any
consistent scoring function.

However, in real world settings, as pointed out by Patton (2015),
forecasts are hardly ever ideal, and the ranking of competing
forecasts might depend on the choice of the scoring function.  This
had already been observed by Murphy (1977), Schervish (1989), and
Merkle and Steyvers (2013), among others, in the important special
case of a binary predictand, where $y = 1$ corresponds to a success
and $y = 0$ to a non-success, so that the mean of the predictive
distribution provides a probability forecast for a success.  As there
is no obvious reason for a consistent scoring function to be preferred
over any other, this raises the question which one of the many
alternatives to use.

Our work is motivated by the quest for guidance in this setting.
Theoretically, the respective key result is that, subject to
unimportant regularity conditions, any function of the form
(\ref{eq:Bregman}) admits a mixture representation of the form
\[
\myS(x,y) = \int_{-\infty}^{+\infty} \myS_\theta(x,y) \: \dd H(\theta), 
\]
where $H$ is a nonnegative measure, and
\begin{eqnarray*}
\myS_\theta(x,y) & = & \left( y - \theta \right)_+ - \left( x - \theta \right)_+ 
                     - \one(x > \theta) \hsp (y - x) \\
                 & = & \begin{cases} 
                       |y - \theta|, & \min(x,y) \leq \theta < \max(x,y), \\
                       0,            & \rm{otherwise}                     
                       \end{cases}
\end{eqnarray*}
for $\theta \in \real$.  (Here and in what follows, we write $(t)_+ =
\max(t,0)$ for the positive part of $t \in \real$ and $\one(A)$ for
the indicator function of the event $A$.)  Thus every scoring function
consistent for the mean can be written as a weighted average over {\em
elementary}\/ or {\em extremal}\/ scores $\myS_\theta$.  As an
important consequence, a point forecast that is preferable in terms of
each extremal score $\myS_\theta$ is preferable in terms of any
consistent scoring function. The elementary scores can be seen as
representing the loss, relative to an oracle, in an investment problem
with cost basis $\theta$ and future revenue $y$; see Section 2.3.

\begin{figure}[t]
	
\centering

\[
\begin{array}{cc}
\toprule
\text{\small Score} & \text{\small Score Difference \rule{0mm}{4mm}} \\
\midrule
\multicolumn{2}{c}{\text{\small Mean Wind Speed, Gneiting et al.~(2006)\rule{0mm}{7mm}}} \\
\includegraphics[height = 4.5cm]{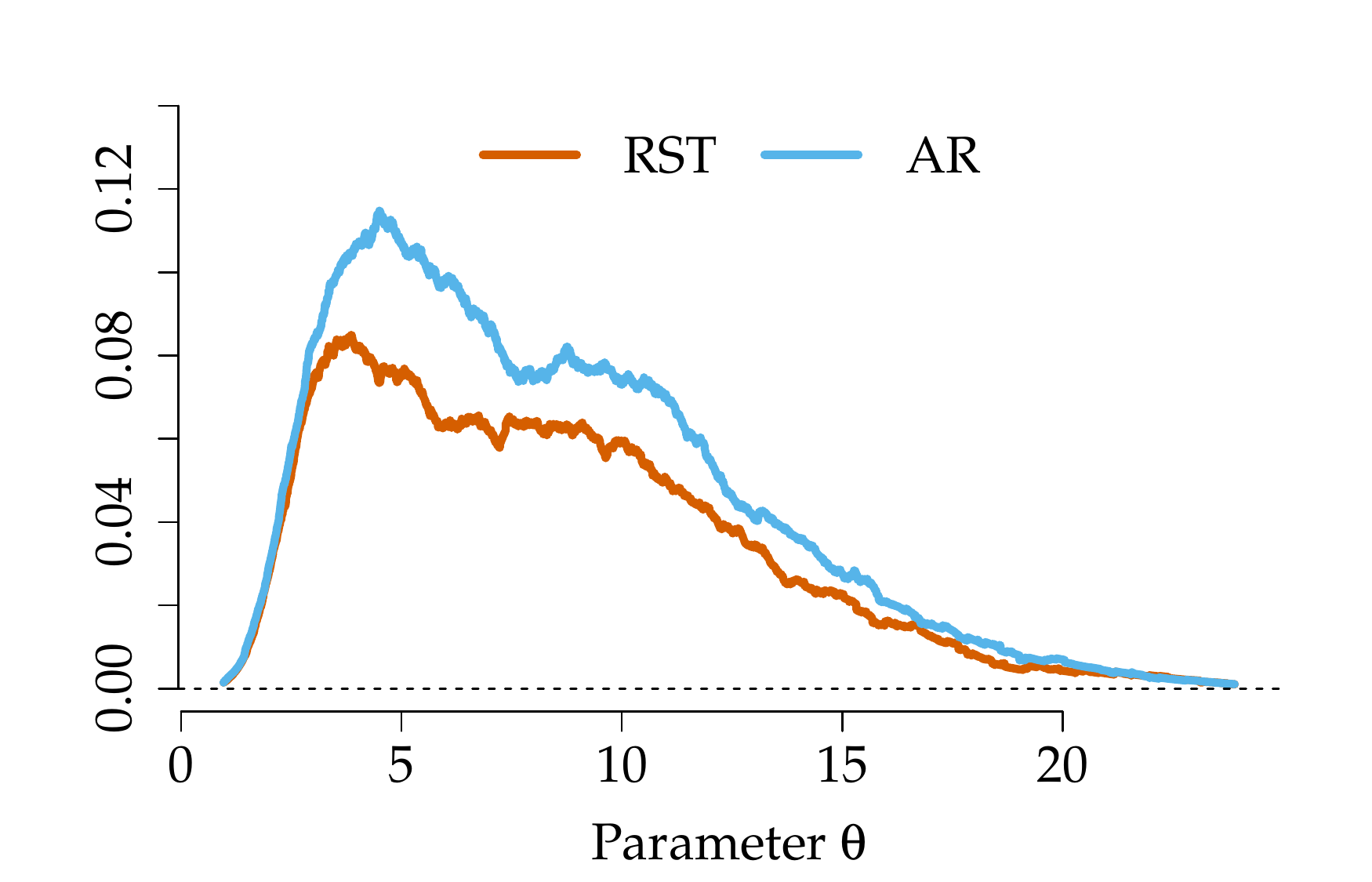} &
\includegraphics[height = 4.5cm]{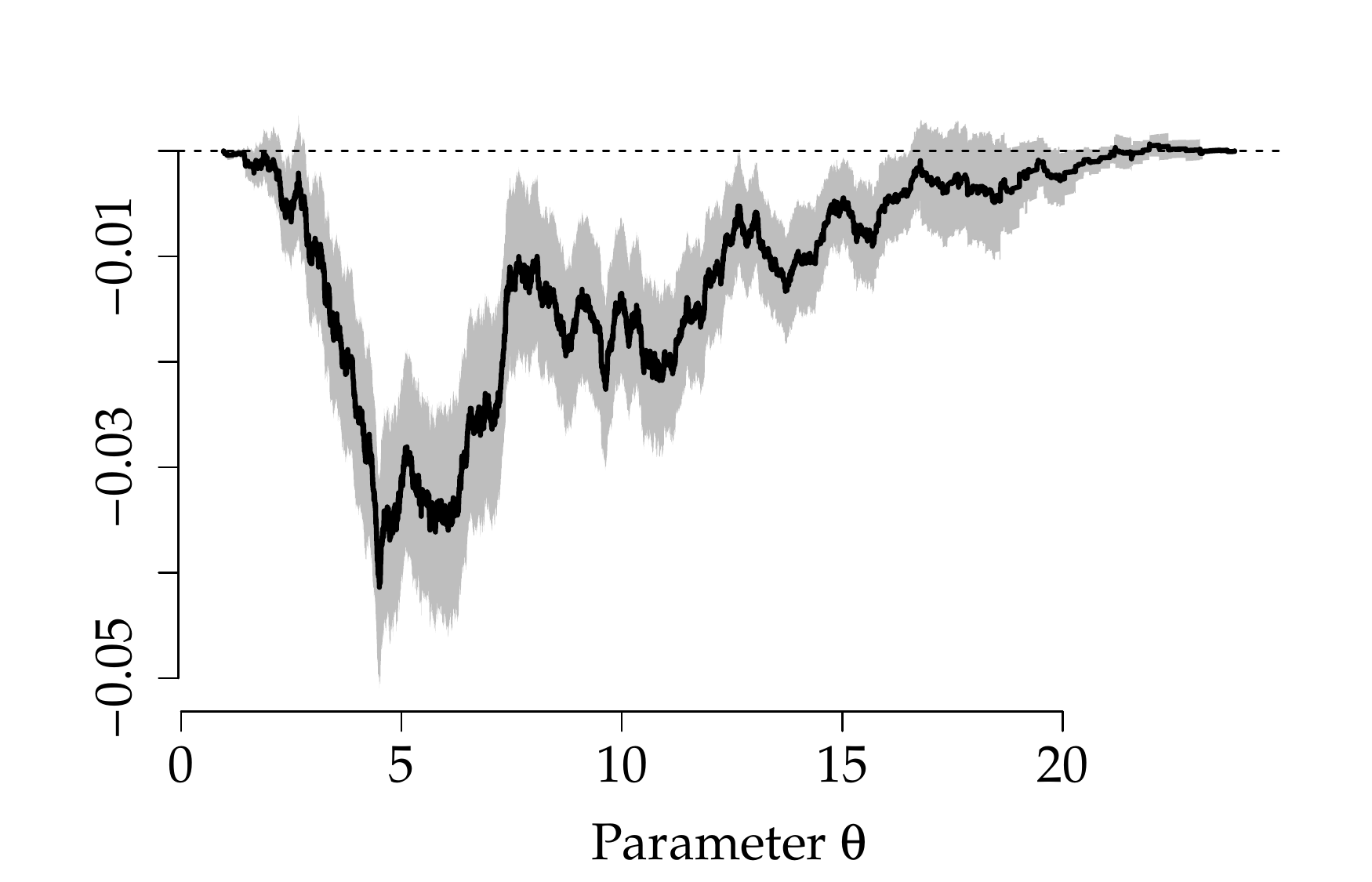} \\ 
\bottomrule
\end{array}
\]

\caption{\small Murphy diagrams for the comparison of point forecasts
  of wind speed at the Stateline wind energy center, using a
  regime-switching space-time (RST) or autoregressive (AR) technique
  (Gneiting et al.~2006).  The functional considered is the mean of
  the respective predictive distribution.  Left: Empirical scores
  $s_j(\theta)$ in (\ref{eq:sc}) versus $\theta$.  Right: Score
  differences along with pointwise 95\% confidence bands.  A negative
  difference means that the RST forecast is preferable.  For details,
  see Sections \ref{subsec:empirical} and \ref{subsec:wind}.
  \label{fig:wind}}

\end{figure}

In empirical settings, point forecasts are compared based on their
average scores.  Specifically, let us consider a sequence of triplets
$(x_{i1}, x_{i2}, y_i)$ for $i = 1, \ldots, n$, where $x_{i1}$ and
$x_{i2}$ are competing point forecasts and $y_i$ is the subsequent
outcome.  We may compare the two forecasts graphically, by plotting
the respective empirical scores,
\begin{equation}  \label{eq:sc} 
s_j(\theta) = \frac{1}{n} \sum_{i=1}^n \myS_\theta(x_{ij},y_i)  
\end{equation} 
for $j = 1$ and $2$, versus $\theta$.  An example of this type of
display, which we term a {\em Murphy diagram}, is shown in Figure
\ref{fig:wind}, where we consider point forecasts of wind speed at a
major wind energy center.

More generally, for both quantiles and expectiles the apparent wealth
of consistent scoring functions can be reduced to a one-dimensional
family of readily interpretable elementary scores, in the sense that
every consistent scoring function can be represented as a mixture from
that family.  The case of the mean or expectation functional, which
includes probability forecasts for binary events as a further special
case, corresponds to the expectile at level $\alpha = 1/2$.

The remainder of the paper is organized as follows.  Section
\ref{sec:theory} is devoted to the key theoretical development, in
which we state and discuss the mixture representations, relate to
Choquet theory and order sensitivity, and provide economic
interpretations of the elementary scores and the associated
functionals.  In particular, we show that expectiles are optimal
decision thresholds in binary investment problems with fixed cost
basis and differential taxation of profits versus losses.  In Section
\ref{sec:rankings}, we apply the mixture representations to study
forecast rankings and propose the aforementioned Murphy diagram for
forecast comparisons.  Illustrations on data examples follow in
Section \ref{sec:examples}, where we revisit meteorological and
economic case studies in the work of Gneiting et al.~(2006), Rudebusch
and Williams (2009), and Patton (2015).  The paper closes with a
discussion in Section \ref{sec:discussion}.  Proofs and computational
details are deferred to Appendices.

\section{Consistent scoring functions for quantiles and expectiles}  \label{sec:theory}

Before focusing on the specific cases of quantiles and expectiles, we
review general background material on the assessment of point
forecasts, with emphasis on consistent scoring functions.

\subsection{Consistent scoring functions}  \label{subsec:consistent}

We first introduce notation and expose conventions.  Let $\cF_0$
denote the class of the probability measures on the Borel-Lebesgue
sets of the real line, $\real$.  For simplicity, we do not distinguish
between a measure $F \in \cF_0$ and the associated cumulative
distribution function (CDF).  We follow standard conventions and
assume CDFs to be right-continuous.  A function $\myS$ defined on a
rectangle $D = D_1 \times D_2 \subseteq \real^2$ is called a {\em
  scoring function}\/ if $\myS(x,y) \geq 0$ for all $(x,y) \in D$ with
$\myS(x,y) = 0$ if $x = y$.  Here, $\myS(x,y)$ is interpreted as the
loss or cost accrued when the point forecast $x$ is issued and the
observation $y$ realizes.  The scoring function is {\em regular}\/ if
it is jointly measurable and left-continuous in its first argument,
$x$, for every $y$.

In point prediction problems, it is rarely evident which functional of
the predictive distribution should be reported.  Guidance can be given
implicitly, by specifying a loss function, or explicitly, by
specifying a functional. The notion of consistency originates in this
setting.

Consider a functional $F \mapsto \myT(F) \subseteq D_1$ on a class
$\cF \subseteq \cF_0$ on which the mapping is well-defined.  Usually,
the functional is single-valued, as in the case of the mean functional
where we take $\cF$ as the class $\cF_1$ of the probability measures
with finite first moment.  More generally, the {\em expectile\/} at
level $\alpha \in (0,1)$ of a probability measure $F \in \cF_1$ is the
unique solution $t$ to the equation 
\[ 
(1-\alpha) \int_{-\infty}^t (t-y) \, \dd F(y) = 
\alpha \int_t^\infty (y-t) \, \dd F(y) \hsp , 
\]
where $\alpha = 1/2$ corresponds to the mean functional (Newey and
Powell 1987).  In the case of quantiles, the functional might be
set-valued.  Specifically, the {\em quantile}\/ functional at level
$\alpha \in (0,1)$ maps a probability measure $F$ to the closed
interval $[q_{\alpha,F}^-, q_{\alpha,F}^+]$, with lower limit
$q_{\alpha,F}^- = \sup \, \{ s : F(s) < \alpha \}$ and upper limit
$q_{\alpha,F}^+ = \sup \, \{ s : \, F(s) \leq \alpha \}$.  The two
limits differ only when the level set $F^{-1}(\alpha)$ contains more
than one point, so typically the functional is {single}-valued.  Any
number between $q_{\alpha,F}^-$ and $q_{\alpha,F}^+$ represents an
$\alpha$-quantile and will be denoted $q_{\alpha,F}$.

The scoring function $\myS$ is {\em consistent}\/ for a functional
$\myT$ relative to the class $\cF$ if {
\begin{equation}  \label{eq:consistent} 
\myE_F \hsp [\myS(t,Y)] \leq \, \myE_F \hsp [\myS(x,Y)] 
\end{equation} 
for all probability measures $F \in \cF$, all $t \in \myT(F)$, and all
point forecasts $x \in D_1$.}  A functional $\myT$ that admits a
strictly consistent scoring function is called {\em elicitable}, and
can then be represented as the solution to an optimization problem, in
that
\[ 
\myT(F) = {\textstyle \arg \min_{\, x}} \, \myE_F \hsp [\myS(x,Y)] \hsp .     
\]
Hence, if the goal is to minimize expected loss, the optimal strategy
is to follow the requested directive in the form of a functional.

In what follows, we restrict attention to the quantile and expectile
functionals.  These are critically important in a gamut of
applications, including quantile and expectile regression in general,
and least squares (i.e., mean) and probit and logit (i.e., binary
probability) regression in particular.

\subsection{Mixture representations} \label{subsec:Choquet}

The classes of the consistent scoring functions for quantiles and
expectiles have been described by Savage (1971), Thomson (1979), and
Gneiting (2011), and we review the respective characterizations in
the setting of the latter paper, where further detail is available.

Up to mild regularity conditions, a scoring function $\myS$ is
consistent for the quantile functional at level $\alpha \in (0,1)$
relative to the class $\cF_0$ if and only if it is of the form
\begin{equation}  \label{eq:qsf}
\myS(x,y) = \left( \one(y < x) -\alpha \right) \left( \hsp g(x) - g(y) \right) \! ,
\end{equation}
where $g$ is non-decreasing.  The most prominent example arises when 
$g(t) = t$, which yields the asymmetric piecewise linear scoring function,
\begin{equation}  \label{eq:apl} 
\myS(x,y) =  \left\{ \begin{array}{ll} 
(1-\alpha) \, (x - y), & y < x, \\ 
\alpha \, (y - x), & y \geq x,
\end{array} \right.
\end{equation} 
that lies at the heart of quantile regression (Koenker and Bassett 1978;
Koenker 2005).  Similarly, a scoring function is consistent for the
expectile at level $\alpha \in (0,1)$ relative to the class $\cF_1$
if and only if it is of the form
\begin{equation}  \label{eq:esf}
\myS(x,y) = |\one(y < x) -\alpha| 
\left( \phi(y) - \phi(x) - \phi'(x) (y-x) \right) \! ,
\end{equation} 
where $\phi$ is convex with subgradient $\phi'$.  The key example
arises when $\phi(t) = t^2$, where
\begin{equation}  \label{eq:ase}  
\myS(x,y) =  \left\{ \begin{array}{ll} 
(1-\alpha) \, (x - y)^2, & y < x. \\ 
\alpha \, (x - y)^2, & y \geq x.
\end{array} \right.
\end{equation} 
This is the loss function used for estimation in expectile regression
(Newey and Powell 1987; Efron 1991), including the ubiquitous case
$\alpha = 1/2$ of ordinary least squares regression.  

In view of the representations (\ref{eq:qsf}) and (\ref{eq:esf}), the
scoring functions that are consistent for quantiles and expectiles are
parameterized by the non-decreasing functions $g$, and the convex
functions $\phi$ with subgradient $\phi'$, respectively.  In general,
neither $g$ nor $\phi$ and $\phi'$ are uniquely determined.  We
therefore select special versions of these functions.  Furthermore, in
the interest of simplicity we generally assume that $D_1 \times D_2 =
\real^2$, adding comments in cases where there are finite boundary
points.  Let $\cI$ denote the class of all left-continuous
non-decreasing real functions, and let $\cC$ denote the class of all
convex real functions $\phi$ with subgradient $\phi' \in \cI$.  This
last condition is satisfied when $\phi'$ is chosen to be the left-hand
derivative of $\phi$, which exists everywhere and is left-continuous
by construction.

In what follows, we use the symbol $\myCQ$ to denote the class of the
scoring functions $\myS$ of the form (\ref{eq:qsf}) where $g \in \cI$.
Similarly, we write $\myCE$ for the class of the scoring functions
$\myS$ of the form (\ref{eq:esf}) where $\phi \in \cC$.  For all
practical purposes, the families $\myCQ$ and $\myCE$ can be identified
with the classes of the regular scoring functions that are consistent
for quantiles and expectiles, respectively.  These classes appear to
be rather large.  However, in either case the apparent multitude can
be reduced to a one-dimensional family of elementary scoring
functions, in the sense that every consistent scoring function admits
a representation as a mixture of elementary elements.

\medskip \noindent {\bf Theorem 1a (quantiles).} 
{\em Any member of the class\/ $\myCQ$ admits a representation of the form
\begin{equation}  \label{eq:mr4q}
\myS(x,y) = \int_{-\infty}^{+\infty} \mySQ(x,y) \, \dd H(\theta)
\qquad (x, y \in \real),  
\end{equation} 
where\/ $H$ is a nonnegative measure and 
\begin{eqnarray}  
\mySQ(x,y) \ceq \left( \one(y < x) - \alpha \right) \, 
\left( \one(\theta < x) - \one(\theta < y) \right) \nonumber \\ \label{eq:qsfx}
\ceq \left\{
\begin{array}{ll} 
1 - \alpha, & y \leq \theta < x, \\  
\alpha,     & x \leq \theta < y, \\ 
0,          & \rm{otherwise}.
\end{array} \right.  
\end{eqnarray} 
The mixing measure\/ $H$ is unique and satisfies\/ $\dd H(\theta) = \dd
\hsp g(\theta)$ for\/ $\theta \in \real$, where\/ $g$ is the nondecreasing
function in the representation\/ (\ref{eq:qsf}).  Furthermore, we have\/
$H(x) - H(y) = \myS(x,y)/(1-\alpha)$ for $x > y$.}

\medskip \noindent {\bf Theorem 1b (expectiles).} 
{\em Any member of the class\/ $\myCE$ admits a representation of the form
\begin{equation}  \label{eq:mr4e}
\myS(x,y) = \int_{-\infty}^{+\infty} \mySE(x,y) \, \dd H(\theta)
\qquad (x, y \in \real),  
\end{equation} 
where\/ $H$ is a nonnegative measure and 
\begin{eqnarray}  
\mySE(x,y) \ceq |\one(y < x) - \alpha|
\left( (y-\theta)_+ - (x-\theta)_+ - (y-x) \, \one(\theta < x) \right) 
\nonumber \\ \label{eq:esfx}
\ceq \left\{
\begin{array}{ll} 
(1 - \alpha) \, |y - \theta|, & y \leq \theta < x, \\  
\alpha \, |y - \theta|,       & x \leq \theta < y, \\ 
0,                            & \rm{otherwise}.
\end{array} \right.  
\end{eqnarray} 
The mixing measure\/ $H$ is unique and satisfies\/ $\dd H(\theta) = \dd
\hsp \phi'(\theta)$ for\/ $\theta \in \real$, where\/ $\phi'$ is the
left-hand derivative of the convex function\/ $\phi$ in the
representation (\ref{eq:esf}).  Furthermore, we have\/ $H(x) - H(y) =
\partial_2 \myS(x,y)/(1-\alpha)$ for\/ $x > y$, where\/ $\partial_2$
denotes the left-hand derivative with respect to the second argument.}

\medskip
Note that the relations in (\ref{eq:mr4q}) and (\ref{eq:mr4e}) hold
pointwise.  In particular, the respective integrals are pointwise
well-defined.  This is because for $(x,y) \in \real^2$ the functions
$\theta \mapsto \mySQ(x,y)$ and $\theta \mapsto \mySE(x,y)$ are
right-continuous, non-negative, and uniformly bounded with bounded
support, and because the non-decreasing functions $g$ and $\phi'$
define non-negative measures $\dd \hsp g$ and $\dd \phi'$ that assign
finite mass to any finite interval.

In the case of quantiles, the asymmetric piecewise linear scoring
function corresponds to the choice $g(t) = t$ in (\ref{eq:qsf}), so
the mixing measure $H$ in the representation (\ref{eq:mr4q}) is the
Lebesgue measure.  The elementary scoring function $\mySQ$ arises when
$g(t) = \one(\theta < t)$, i.e., when $H$ is a one-point measure in
$\theta$.

In the case of expectiles, the mixing measure for the asymmetric
squared error scoring function is twice the Lebesgue measure.  The
choice $\alpha = 1/2$ recovers the mean or expectation functional, for
which existing parametric subfamilies emerge as special cases of our
mixture representation.  Patton's (2015) exponential Bregman family,
\[
\myS_a(x,y) = 
\frac{1}{a^2} \left( \exp(ay) - \exp(ax) \right) 
- \frac{1}{a} \exp(ax) \hsp (y-x)
\qquad (a \neq 0),
\]
which nests the squared error loss in the limit as $a \to 0$,
corresponds to the choice $\phi(t) = a^{-2} \exp(at)$ in
(\ref{eq:esf}).  The mixing measure $H$ in the representation
(\ref{eq:mr4e}) then has Lebesgue density $h(\theta) = \exp(a\theta)$
for $\theta \in \real$.  For Patton's (2011) family
\[
\myS_b(x,y) = \left\{ \begin{array}{ll}
\displaystyle 
\frac{y^b - x^b}{b \hsp (b-1)} - \frac{x^{b-1}}{b-1} \, (y-x), & b \notin \{ 0, 1 \}, \\ 
\displaystyle \rule{0mm}{6.5mm}
\frac{y}{x} - \log \frac{y}{x} - 1,                           & b = 0, \\ 
\displaystyle \rule{0mm}{6.5mm}
y \log \frac{y}{x} - (y-x),                                   & b = 1,
\end{array} \right.
\]
of homogeneous scoring functions on the positive half line the mixing
measure has Lebesgue density $h(\theta) = \theta^{\, b-2} \hsp
\one(\theta > 0)$, remarkably with no case distinction being required.
The elementary scoring function $\mySE$ emerges when $\phi(t) = (t -
\theta)_+$ in (\ref{eq:esf}); here the mixing measure in
(\ref{eq:mr4e}) is a one-point measure in $\theta$.

From a theoretical perspective, a natural question is whether the
mixture representations (\ref{eq:mr4q}) and (\ref{eq:mr4e}) can be
considered {\em Choquet representations}\/ in the sense of functional
analysis (Phelps 2001).  Recall that a member $\myS$ of a convex class
$\cS$ is an {\em extreme point}\/ of $\cS$ if it cannot be written as
an average of two other members, i.e., if $\myS = (\myS_1+ \myS_2)/2$
with $\myS_1, \myS_2 \in \cS$ implies $\myS_1 = \myS_2 = \myS$.  Our
mixture representations qualify as Choquet representations if the
elementary scores $\mySQ$ and $\mySE$ form extreme points of the
underlying classes of scoring functions.  This cannot possibly be true
for our classes $\myCQ$ and $\myCE$ because they are invariant under
dilations, hence admit trivial average representations built with
multiples of one and the same scoring function. Therefore, the
families $\myCQ$ and $\myCE$ need to be restricted suitably.
Specifically, let the class $\cI_1$ consist of all functions $g \in
\cI$ such that $\lim_{x \to -\infty} g(x) = 0$ and $\lim_{x \to
+\infty} g(x) = 1$.  Similarly, let $\cC_1$ denote the family of all
$\phi \in \cC$ such that $\phi(0) = 0$ and $\phi' \in \cI_1$.  These
classes are convex, and so are the associated subclasses of the
families $\myCQ$ and $\myCE$, which we denote by $\cS_{\alpha,1}^{\rm
Q}$ and $\cS_{\alpha,1}^{\rm E}$, respectively.  The elementary scores
$\mySQ$ and $\mySE$ evidently are members of these restricted
families.

\medskip \noindent {\bf Proposition 1a (quantiles).} 
{\em For every\/ $\alpha \in (0,1)$ and\/ $\theta \in \real$, the
scoring function\/ $\mySQ$ is an extreme point of the class\/
$\cS_{\alpha,1}^{\rm Q}$.}

\medskip \noindent {\bf Proposition 1b (expectiles).} 
{\em For every\/ $\alpha \in (0,1)$ and\/ $\theta \in \real$, the 
scoring function\/ $\mySE$ is an extreme point of the class\/ 
$\cS_{\alpha,1}^{\rm E}$.}

\medskip 
We thus have furnished Choquet representations for subclasses of
the consistent scoring functions for quantiles and expectiles.  In the
extant literature, such Choquet representations have been known in the
binary case only, where $y = 1$ corresponds to a
success and $y = 0$ to a non-success, so that the mean, $p \in
[0,1]$, of the predictive distribution provides a probability forecast
for a success.  In this setting, the Savage representation
(\ref{eq:esf}) for the members of the respective class
$\cS_{1/2,1}^{\rm E}$ reduces to
\[
\myS( \hsp p,0) = 
\frac{1}{2} \left( \hsp p \hsp \phi'( \hsp p) - \phi( \hsp p) \right) \! , \quad 
\myS( \hsp p,1) = 
\frac{1}{2} \left( \phi(1) - \phi( \hsp p) - (1-p) \hsp \phi'( \hsp p) \right) \! . 
\]
The mixture representation (\ref{eq:mr4e}) can then be written as 
\begin{equation}  \label{eq:mr4p}
\myS( \hsp p,y) = \int_0^1 \mySB( \hsp p,y) \, \dd H(\theta), 
\end{equation} 
where $H$ is a nonnegative measure and 
\begin{equation}  \label{eq:psfx}  
\mySB( \hsp p,y) 
= 2 \hsp \myS^{\rm E}_{1/2,\theta}( \hsp p,y) 
= \left\{ \begin{array}{ll} 
  \theta, & y = 0, \; p > \theta, \\  
  1 - \theta,     & y = 1, \; p \leq \theta, \\ 
  0,          & \rm{otherwise}. 
  \end{array} \right.   
\end{equation}

The parameter $\theta \in (0,1)$ can be interpreted as the cost-loss
ratio in the classical simple cost-loss decision model (Richardson
2012).  Up to unimportant conventions regarding coding, scaling, and
gain-loss orientation, this recovers the well known mixture
representation of the proper scoring rules for probability forecasts
of binary events (Shuford, Albert, and Massengill 1966; Schervish
1989).  Different choices of the mixing measure yield the standard
examples of scoring rules in this case; see Buja et al.~(2005) and
Table 1 in Gneiting and Raftery (2007).  The widely used Brier score,
\begin{equation}  \label{eq:Brier}  
\myS( \hsp p,0) = p^2, \quad 
\myS( \hsp p,1) = (1-p)^2,  
\end{equation} 
arises when $H$ is twice the Lebesgue measure. 

We close the section by noting a fundamental connection between the
extremal scoring rules for quantiles, expectiles, and probabilities in
(\ref{eq:qsfx}), (\ref{eq:esfx}), and (\ref{eq:psfx}), respectively.
Specifically, given any predictive CDF, $F$, and outcome, $y \in \real$,
\begin{equation}  \label{eq:qep}  
\mySQ \! \left( q_{\alpha,F}^-, y \right) 
= 2 \hsp \myS^{\rm E}_{1/2,1-\alpha} \! \left( 1-F(\theta), \one(y > \theta) \right)   
\end{equation}
for every $\alpha \in (0,1)$ and $\theta \in \real$.  We will revisit this 
relation repeatedly.

\subsection{Economic interpretation}  \label{subsec:econ}

Our results in the previous section give rise to natural economic
interpretations of the extremal scoring functions $\mySQ$ and $\mySE$,
along with the quantile and expectile functionals themselves.  In
either case, the interpretation relates to a binary betting or
investment decision with random outcome, $y$.

In the case of the extremal quantile scoring function $\mySQ$ in
(\ref{eq:qsfx}), the payoff takes on only two possible values,
relating to a bet on whether or not the outcome $y$ will exceed the
threshold $\theta$.  Specifically, consider the following payoff
scheme, which is realized in spread betting in prediction markets
(Wolfers and Zitzewitz 2008):

\begin{itemize} 

\item If Quinn refrains from betting, his payoff
  will be zero, independently of the outcome $y$.

\item If Quinn enters the bet and $y \leq \theta$
  realizes, he loses his wager, $\rho_L > 0$.

\item If Quinn enters the bet and $y > \theta$ realizes, his winnings
  are $\rho_G > \rho_L$, for a gain of $\rho_G - \rho_L$.

\end{itemize} 

How should Quinn act under this payoff scheme?  If Quinn does not
enter the bet, his actual and expected payoffs equal zero.  If he does
enter, his expected payoff is
\[
- \rho_L \int_{-\infty}^\theta \dd F(y) 
+ \left( \hsp \rho_G - \rho_L \right) \int_\theta^\infty \dd F(y),  
\] 
where $F$ is Quinn's predictive CDF for the future outcome, $y$, which
for simplicity we assume to be strictly increasing.  This
expression is strictly positive if and only if $q_{\alpha,F} >
\theta$, where
\begin{equation}  \label{eq:alphaq}
\alpha = \frac{\rho_G - \rho_L}{\rho_G} \in (0,1). 
\end{equation} 
Hence, Quinn's optimal decision rule is determined by the
$\alpha$-quantile of $F$, in that he enters the bet if and only if
$q_{\alpha,F} > \theta$.  Motivated by the specific format of the
optimal decision or Bayes rule, the top left matrix in Table
\ref{tab:payoffs} summarizes the payoff from just any strategy of the
form {\em enter the bet if and only if}\/ $x > \theta$.

It remains to draw the connection to the extremal scoring function
$\mySQ$.  To this end, we shift attention from positively oriented
payoffs to negatively oriented regrets, which we define as the
difference between the payoff for an oracle and Quinn's payoff.  Here
the term oracle refers to a (hypothetical) omniscient bettor who
enters the bet if and only if $y > \theta$ realizes, which would yield
an ideal payoff $\rho_G - \rho_L$ if $y > \theta$, and zero otherwise.
If Quinn uses some decision threshold $x$, his regret equals the
extremal score $\mySQ(x,y)$ except for an irrelevant multiplicative
factor.  This is illustrated in the bottom left matrix in the table
and corresponds to the classical, simple cost-loss decision model
(Richardson 2012).  In decision theoretic terms, the distinction
between payoff and regret is inessential, because the difference
depends on the outcome, $y$, only.  In either case, the optimal
strategy is to choose the decision threshold $x = q_{\alpha,F}$.

\begin{table}[t]

\centering

\small

\begin{tabular}{cc}
\toprule
Quantiles              & Expectiles \\ 
\midrule 
Monetary Payoff & Monetary Payoff \rule{0mm}{6.5mm} \\ [2mm]
\begin{tabular}{c|cc} 
                & $y \leq \theta$ & $y > \theta$ \\ [1mm]
\hline
$x \leq \theta$ & 0               & 0 \rule{0mm}{4.5mm} \\
$x > \theta$    & $-\rho_L$       & $\rho_G - \rho_L$ \rule{0mm}{3.5mm} \\
\end{tabular} 
& 
\begin{tabular}{c|cc} 
                & $y \leq \theta$ & $y > \theta$ \\ [1mm] 
\hline
$x \leq \theta$ & 0               & 0 \rule{0mm}{4.5mm} \\
$x > \theta$    & $-(1-\kappa_L)(\theta-y)$ & $(1-\kappa_G)(y-\theta)$ \rule{0mm}{3.5mm} \\
\end{tabular} 
\\
Score (Regret)  & Score (Regret) \rule{0mm}{7.5mm} \\ [2mm]
\begin{tabular}{c|cc} 
                & $y \leq \theta$ & $y > \theta$ \\ [1mm]
\hline
$x \leq \theta$ & 0               & $\rho_G - \rho_L$ \rule {0mm}{4.5mm} \\
$x > \theta$    & $\rho_L$        & 0 \rule{0mm}{3.5mm} \\
\end{tabular} 
& 
\begin{tabular}{c|cc} 
                & $y \leq \theta$ & $y > \theta$ \\ [1mm] 
\hline
$x \leq \theta$ & 0               & $(1-\kappa_G)(y-\theta)$ \rule {0mm}{4.5mm} \\
$x > \theta$    & $(1-\kappa_L)(\theta-y)$ & 0 \rule{0mm}{3.5mm} \\
\end{tabular} 
\\ \\
\bottomrule
\end{tabular} 

\caption{\small Overview of payoff structures for decision rules of
  the form {\em enter the bet/invest if and only if\/} $x > \theta$.
  Monetary payoffs are positively oriented, whereas scores are
  negatively oriented regrets relative to an oracle.  In the left
  column, the regret equals the extremal score $\mySQ(x,y)$, where
  $\alpha = (\rho_G - \rho_L) / \rho_G$, up to a multiplicative
  factor.  In the right column, the regret is $\mySE(x,y)$, where
  $\alpha = (1 - \kappa_G) / (2 - \kappa_G - \kappa_L)$, again up to a
  multiplicative factor.  \label{tab:payoffs}}

\end{table}

In the case of the extremal expectile scoring function $\mySE$ in
(\ref{eq:esf}), the payoff is real-valued.  Specifically, suppose that
Eve considers investing a fixed amount $\theta$ into a start-up
company, in exchange for an unknown, future amount $y$ of the
company's profits or losses.  The payoff structure then is as follows:

\begin{itemize}

\item If Eve refrains from the deal, her payoff will
  be zero, independently of the outcome $y$.

\item If Eve invests and $y \leq \theta$ realizes, her payoff is
  negative, at $- \hsp (1-\kappa_L) \hsp (\theta - y)$. Here, $\theta
  - y$ is the sheer monetary loss, and the factor $1 - \kappa_L$
  accounts for Eve's reduction in income tax, with $\kappa_L \in
  [0,1)$ representing the deduction rate.\footnote{In financial terms,
  the loss acts as a tax shield. The linear functional form assumed
  here is not unrealistic, even though it is simpler than many
  real-world tax schemes, where nonlinearities may arise from tax
  exemptions, progression, etc.}
 
\item If Eve invests and $y > \theta$ realizes, her payoff is
  positive, at $(1-\kappa_ G) \hsp (y-\theta)$, where $\kappa_G \in
  [0, 1)$ denotes the tax rate that applies to her profits.

\end{itemize}

How should Eve act under this payoff scheme?  If Eve does not enter
the deal, her actual and expected payoffs vanish.  In case she
invests, the expected payoff is
\[
- \left( 1 - \kappa_L \right) \int_{-\infty}^\theta (\theta - y) \, \dd F(y) 
+ \left( 1 - \kappa_G \right) \int_\theta^\infty (y - \theta) \, \dd F(y).  
\] 
This expression is strictly positive if and only if the expectile at level 
\begin{equation}  \label{eq:alphae}
\alpha = \frac{1 - \kappa_G}{2 - \kappa_G -\kappa_L} \in (0,1)  
\end{equation} 
of Eve's predictive CDF, $F$, exceeds $\theta$.  In analogy to the
quantile case, the top right matrix in Table \ref{tab:payoffs}
represents Eve's payoff from just any strategy of the form {\em invest
if and only if}\/ $x > \theta$.

To relate to the extremal scoring function $\mySE$, we again shift
attention to regrets relative to an omniscient investor or oracle who
enters the deal if and only if $y > \theta$ occurs, which would yield
the ideal payoff $(1 - \kappa_G) (y - \theta)_+$.  As seen in the
table, if Eve uses the threshold $x$ to determine whether or not to
invest, the regret equals the extremal score $\mySE(x,y)$, up to a
multiplicative factor.\footnote{The elementary score $\mySB$ for
probability forecasts of a binary event in (\ref{eq:psfx}) is obtained
when $\kappa_G = \kappa_L = 0$ and $y \in \{ 0, 1 \}$.  The parameter
$\theta \in (0,1)$ can then be interpreted as a cost-loss ratio.}

Therefore, expectiles can be interpreted as optimal decision
thresholds in investment problems with fixed costs and differential
tax rates for profits versus losses.  The mean arises in the special
case when $\alpha = 1/2$ in (\ref{eq:alphae}).  It corresponds to
situations in which losses are fully tax deductible ($\kappa_G =
\kappa_L$) and nests situations without taxes ($\kappa_G = \kappa_L =
0$).  Tough taxation settings where $\kappa_L < \kappa_G$ shift Eve's
incentives toward not entering the deal and correspond to expectiles
at levels $\alpha < 1/2$.  For example, if losses cannot be deducted
at all $(\kappa_L = 0)$, whereas profits are taxed at a rate of
$\kappa_G = 1/2$, Eve will invest only if the expectile at level
$\alpha = 1/3$ of her predictive CDF, $F$, exceeds the deal's fixed
costs, $\theta$.  Note that we permit the case $\theta < 0$, which may
reflect subsidies or tax credits, say.

The above interpretation of expectiles as optimal thresholds in
decision problems attaches an economic meaning to this class of
functionals, which thus far seems to have been missing; e.g., Schulze
Waltrup et al.~(2014, p.~2) note that ``expectiles lack an intuitive
interpretation''.  The foregoing may also bear on the debate about the
revision of the Basel protocol for banking regulation, which involves
contention about the choice of the functional of in-house risk
distributions that banks are supposed to report to regulators
(Embrechts et al.~2014).  Recently, expectiles have been put forth as
potential candidates, as it has been proved that they are the only
elicitable law-invariant coherent risk measures (Delbaen et al.~2014;
Ziegel 2014; Bellini and Bignozzi 2015).

\subsection{Order sensitivity}  \label{subsec:ordsen}

The extremal scoring functions $\mySQ$ and $\mySE$ are not only
consistent for their respective functional, they in fact enjoy the
stronger property of order sensitivity.  Generally, a scoring function
$\myS$ is {\em order sensitive}\/ for the 
functional $F \mapsto T(F)$ relative to the class $\cF$ if, for all 
$F \in \cF$, all $t \in \myT(F)$, and all $x_1, x_2 \in \real$,
\[
x_2 \leq x_1 \leq t \; \Longrightarrow \;
\myE_F \hsp [\myS(x_2,Y)] \geq \myE_F \hsp [\myS(x_1,Y)],  
\]
and 
\[
t \leq x_1 \leq x_2 \; \Longrightarrow \; 
\myE_F \hsp [\myS(x_1,Y)] \leq \myE_F \hsp [\myS(x_2,Y)].
\]
The order sensitivity is {\em strict}\/ if the above continues to hold
when the inequalities involving $x_1$ and $x_2$ are strict.  As before, we
denote the class of the Borel probability measures on $\real$ by
$\cF_0$, and we write $\cF_1$ for the subclass of the probability
measures with finite first moment.

\medskip \noindent {\bf Proposition 2a (quantiles).}  
{\em For every\/ $\alpha \in (0,1)$ and\/ $\theta \in \real$, the
extremal scoring function\/ $\mySQ$ is order sensitive for the\/
$\alpha$-quantile functional relative to\/ $\cF_0$.}

\medskip \noindent {\bf Proposition 2b (expectiles).}
{\em For every\/ $\alpha \in (0,1)$ and\/ $\theta \in \real$, the
extremal scoring function\/ $\mySE$ is order sensitive for the\/
$\alpha$-expectile functional relative to\/ $\cF_1$.}
  
\medskip 
Owing to the mixture representations (\ref{eq:mr4q}) and
(\ref{eq:mr4e}), the order sensitivity of the extremal scoring
functions transfers to all regular consistent scoring functions.
Strict order sensitivity applies if the function $g$ in the
representation (\ref{eq:qsf}) and the derivative $\phi'$ in the
representation (\ref{eq:esf}), respectively, are strictly increasing,
relative to subclasses of probability measures with suitable moment
constraints.  Closely related results have recently been obtained in
studies of elicitability (Steinwart et al.~2014; Ziegel 2014; Bellini
and Bignozzi 2015).  In this strand of literature, the ambitious goal
of characterizing all elicitable functionals necessitates regularity
conditions that are not satisfied by our discontinuous, compactly
supported extremal scoring functions.

\section{Forecast rankings}  \label{sec:rankings}

In this section, we turn to the task of comparing and ranking
forecasts.  Before {applying} our mixture representations to this
problem, we introduce the prediction space setting of Gneiting and
Ranjan (2013) and define notions of forecast dominance.

\subsection{Prediction spaces}  \label{subsec:space} 

A {\em prediction space}\/ is a probability space tailored to the
study of forecasting problems.  Following the seminal work of Murphy
and Winkler (1987), the prediction space setting of Gneiting and
Ranjan (2013) considers the joint distribution of forecasts and
observations.  We first focus on probabilistic forecasts, $F$, which
we identify with the associated cumulative distribution functions
(CDFs) for the real-valued outcome, $Y$.  The elements of the
respective sample space $\Omega$ can be identified with tuples of
the form
\begin{equation}  \label{eq:Fspace} 
\left( F_1, \ldots, F_k, Y \right), 
\end{equation} 
where the predictive distributions $F_1, \ldots, F_k$ utilize
information sets $\cA_1, \ldots, \cA_k \subseteq \cA$, respectively,
with $\cA$ being a sigma field on the sample space $\Omega$.  In
measure theoretic language, the information sets correspond to sub
sigma fields, and $F_j$ is a CDF-valued random quantity measurable
with respect to $\cA_j$.  The joint distribution of the quantities in
(\ref{eq:Fspace}) is encoded by a probability measure $\myQ$ on
$(\Omega, \cA)$.  In this setting, a predictive distribution $F_j$ is
{\em ideal}\/ relative to $\cA_j$ if it corresponds to the conditional
distribution of the outcome $Y$ under $\myQ$ given $\cA_j$.

\begin{table}[t]

\centering 

\footnotesize

\begin{tabular}{llccc}
\toprule
Forecaster & Predictive Distribution & $\alpha$-Quantile & Mean & Prob$(Y > y)$ \rule{0mm}{3mm}\\
\midrule
Perfect        & $\cN(\mu,1)$ & $\mu + z_\alpha$ & $\mu$ & $1- \Phi(y - \mu)$ \rule{0mm}{4mm} \\
Climatological & $\cN(0,2)$   & $\sqrt{2} \hsp z_\alpha$ & 0 & $1 - \Phi(\frac{y}{\sqrt{2}})$ \rule{0mm}{5mm} \\ 
Unfocused      & $\frac{1}{2} (\cN(\mu,1) + \cN(\mu+\tau,1))$ 
               & $\mu + z_{\alpha,\tau}$ & $\mu + \frac{\tau}{2}$ & $1-\Phi_\tau(y - \mu)$ \rule{0mm}{5mm} \\
Sign-reversed  & $\cN(-\mu,1)$ & $- \mu + z_\alpha$ & $- \mu$ & $1 - \Phi(y + \mu)$ \rule{0mm}{5mm} \\
\bottomrule
\end{tabular}

\caption{\small An example of a prediction space with four competing
  forecasters.  The outcome is generated as $Y \, | \, \mu \sim
  \cN(0,1)$, where $\mu \sim \cN(0,1)$.  The random variable $\tau$
  attains the values $-2$ and 2 with probability $1/2$, independently
  of $\mu$ and $Y$.  For $\alpha \in (0,1)$ and $\tau \in \{ -2, 2
  \}$, we let $z_\alpha = \Phi^{-1}(\alpha)$, $\Phi_\tau(x) = (\Phi(x)
  + \Phi(x-\tau))/2$, and $z_{\alpha,\tau} = \Phi_\tau^{-1}(\alpha)$,
  where $\Phi$ denotes the CDF of the standard normal
  distribution.  \label{tab:sim}}

\end{table}

In a nutshell, a prediction space specifies the joint distribution of
tuples of the form (\ref{eq:Fspace}).  To give an example, Table
\ref{tab:sim} revisits a scenario studied by Gneiting et al.~(2007)
and Gneiting and Ranjan (2013).\footnote{The only difference is that
we let the random variable $\tau$ attain the values $-2$ and 2, rather
than the values $-1$ and 1 as in Gneiting et al.~(2007) and Gneiting
and Ranjan (2013).}  Here, the outcome is generated as $Y \, | \, \mu
\sim \cN(0,1)$ where $\mu \sim \cN(0,1)$.  The perfect forecaster is
ideal relative to the sigma field generated by the random variable
$\mu$.  The unfocused and sign-reversed forecasters also have
knowledge of $\mu$, but fail to be ideal.  The climatological
forecaster, issuing the unconditional distribution of the outcome $Y$
as predictive distribution, is ideal relative to the uninformative
sigma field generated by the empty set.

Any predictive distribution, $F$, can be reduced to a point forecast
by extracting the sought functional, $\myT(F)$.  In what follows, we
focus on quantiles, the mean or expectation functional, and
probability forecasts of the binary event that the outcome exceeds a
threshold value.  The respective point forecasts for the perfect,
climatological, unfocused, and sign-reversed forecaster are shown in
Table \ref{tab:sim}.

In practice, point forecasts might be an end to themselves, i.e., they
might have been issued without there being an underlying predictive
distribution.  To accommodate such cases, we define a {\em point
  prediction space}\/ to be a probability space $(\Omega, \cA, \myQ)$,
where the elements of the sample space $\Omega$ can be identified with
tuples of the form
\begin{equation}  \label{eq:Xspace} 
\left( X_1, \ldots, X_l, Y \right), 
\end{equation} 
where the random variables $X_1, \ldots, X_l$ represent point
forecasts and utilize information sets $\cA_1, \ldots, \cA_l \subseteq
\cA$, respectively.\footnote{For simplicity, we let $X_1, \ldots, X_l$
  be single-valued.  Extensions to set-valued random
  quantities, as might occur in the case of quantiles, are
  straightforward.}  The joint distribution of the point forecasts and
the observation in (\ref{eq:Xspace}) is specified by the probability
measure $\myQ$.  Similarly, it is sometimes useful to consider a {\em
  mixed prediction space}, by specifying the joint distribution $\myQ$
of tuples of the form
\begin{equation}  \label{eq:FXspace} 
\left( F_1, \ldots, F_k, X_1, \ldots, X_l, Y \right), 
\end{equation} 
where $F_1, \ldots, F_k$ represent CDF-valued random quantities, and
$X_1, \ldots, X_l$ represent point forecasts.

\subsection{Notions of forecast dominance}  \label{subsec:dominance} 

We now define notions of forecast dominance, starting with
probabilistic forecasts that take the form of predictive CDFs, and
then turning to point forecasts. In the former setting, a scoring rule
$\hat{\myS}(F,y)$ is a suitably measurable function that assigns a loss
or penalty when we issue the predictive distribution $F$ and $y$
realizes.  A scoring rule $\hat{\myS}$ is {\em proper}\/ if
\begin{equation}  \label{eq:proper}  
\myE_G \, \hat{\myS}(G,Y) \leq \myE_G \, \hat{\myS}(F,Y)
\end{equation} 
for all probability measures $F$ and $G$ in its domain of definition
(Gneiting and Raftery 2007).  Proper scoring rules therefore encourage
honest and careful assessments.  As is well known, a scoring function
$\myS$ that is consistent for a single-valued functional $\myT$ relative to a
class $\cF$ induces a proper scoring rule, by defining
$\hat{\myS}(F,y) = \myS(\myT(F),y)$ for $F \in \cF$ and $y \in \real$.

\medskip \noindent {\bf Definition 1 (predictive CDFs).}  
{\em Let\/ $F_1$ and\/ $F_2$ be probabilistic forecasts, and let\/ $Y$
be the outcome, in a prediction space.  Then\/ $F_1$ dominates\/ $F_2$
relative to a class\/ $\cP$ of proper scoring rules if\/ $\myE_\myQ \,
\hat{\myS}(F_1,Y) \leq \myE_\myQ \, \hat{\myS}(F_2,Y)$ for every\/
$\hat{\myS} \in \cP$.}

\medskip
We now turn to quantiles and expectiles and the respective families
$\myCQ$ and $\myCE$ of the regular consistent scoring functions for
these functionals.  

\medskip \noindent {\bf Definition 2a (quantiles).}
{\em Let\/ $X_1$ and\/ $X_2$ be point forecasts, and let\/ $Y$ be the
  outcome, in a point prediction space.  Then\/ $X_1$ dominates\/
  $X_2$ as an $\alpha$-quantile forecast if\/ $\myE_\myQ \,
  \myS(X_1,Y) \leq \myE_\myQ \, \myS(X_2,Y)$ for every scoring
  function\/ $\myS \in \myCQ$.}

\medskip \noindent {\bf Definition 2b (expectiles).}
{\em Let\/ $X_1$ and\/ $X_2$ be point forecasts, and let\/ $Y$ be the
  outcome, in a point prediction space.  Then\/ $X_1$ dominates\/
  $X_2$ as an $\alpha$-expectile forecast if\/ $\myE_\myQ \,
  \myS(X_1,Y) \leq \myE_\myQ \, \myS(X_2,Y)$ for every scoring
  function\/ $\myS \in \myCE$.}

\medskip
It is important to note that the expectations in the definitions are
taken with respect to the joint distribution of the probabilistic
forecasts and the outcome.  The notions provide partial orderings for
the predictive distributions $F_1, \ldots, F_k$ in (\ref{eq:Fspace})
and $X_1, \ldots, X_k$ in (\ref{eq:Xspace}), respectively.\footnote{In
  the special case of probability forecasts of a binary event, related
  notions of sufficiency and dominance have been studied by DeGroot
  and Fienberg (1983), Vardeman and Meeden (1983), Schervish (1989),
  Kr\"amer (2005), and Br\"ocker (2009).}  Essentially, a
probabilistic forecast that dominates another is preferable, or at
least not inferior, in any type of decision that involves the
respective predictive distributions.  In the case of quantiles or
expectiles, a point forecast that dominates another is preferable, or
at least not inferior, in any type of decision problem that depends on
the respective predictive distributions via the considered functional
only. Adaptations to functionals other than quantiles or expectiles
are straightforward.
  
Under which conditions does a forecast dominate another?  Holzmann and
Eulert (2014) recently showed that if two predictive distributions are
ideal, then the one with the richer information set dominates the
other.  Furthermore, the result carries over to ideal forecasters'
induced point predictions, including but not limited to the cases of
quantiles and expectiles that we consider here.  To give an example in
the setting of Table \ref{tab:sim}, the perfect and the climatological
forecasters are ideal relative to the sigma fields generated by $\mu$,
and generated by the empty set, respectively.  Therefore, the perfect
forecaster dominates the climatological forecaster, in any of the
above senses.

Tsyplakov (2014) went on to show that if a predictive distribution is
ideal relative to a certain information set, then it dominates any
predictive distribution that is measurable with respect to the
information set.  Again, the result carries over to the induced point
forecasts.  In the setting of Table \ref{tab:sim}, the perfect
forecaster is ideal relative to the sigma field generated by the
random variables $\mu$ and $\tau$.  The climatological, unfocused, and
sign-reversed forecasters are measurable with respect to this sigma
field, and so they are dominated by the perfect forecaster, in any of
the above senses.

In the practice of forecasting, predictive distributions are hardly
ever ideal, and information sets may not be nested, as emphasized by
Patton (2015). Therefore, the above theoretical results are not
readily applicable, and distinct soring rules, or distinct consistent
scoring functions, may yield distinct forecast rankings, as in
empirical examples given by Schervish (1989), Merkle and Steyvers
(2013), and Patton (2015), among others.  Furthermore, in general it
is not feasible to check the validity of the expectation inequalities
in Definitions 1, 2a, and 2b for {\em any}\/ proper scoring rule
$\hat{\myS} \in \cP$, or consistent scoring function $\myS \in \myCQ$,
or $\myS \in \myCE$, respectively.

Fortunately, in the case of quantile and expectile forecasts, the
mixture representations in Theorems 1a and 1b reduce checks for
dominance to the respective one-dimensional families of elementary
scoring functions.

\medskip \noindent {\bf Corollary 1a (quantiles).}  
{\em In a point prediction space, $X_1$ dominates\/ $X_2$ as an\/
  $\alpha$-quantile forecast if\/ $\myE_\myQ \, \mySQ(X_1,Y) \leq
  \myE_\myQ \, \mySQ(X_2,Y)$ for every\/ $\theta \in \real$.}

\medskip \noindent {\bf Corollary 1b (expectiles).}  
{\em In a point prediction space, $X_1$ dominates\/ $X_2$ as an\/
  $\alpha$-expectile forecast if\/ $\myE_\myQ \, \mySE(X_1,Y) \leq
  \myE_\myQ \, \mySE(X_2,Y)$ for every\/ $\theta \in \real$.}

\medskip
The reduction to a one-dimensional problem suggests graphical
comparisons via Murphy diagrams.  Before we discuss this tool, we note
that order sensitivity can sometimes be invoked to prove dominance.
For example, consider the mixed prediction space setting
(\ref{eq:FXspace}) with $k = 2$ and $l = 1$.  Suppose that the
CDF-valued random quantity $F$ is ideal relative to the sigma field
$\cA$, and let $q_{\alpha, F}$ denote its $\alpha$-quantile.  Suppose
furthermore that $X_1$ and $X_2$ are measurable with respect to $\cA$.
By Corollary 1a in concert with Proposition 1a and a conditioning
argument, $X_1$ dominates $X_2$ as an $\alpha$-quantile forecast if 
with probability one either 
\[
X_2 \leq X_1 \leq q_{\alpha,F} 
\qquad \mbox{or} \qquad 
q_{\alpha,F} \leq X_1 \leq X_2
\]
holds true.  An analogous argument applies in the case of the
$\alpha$-expectile.

In the scenario of Table \ref{tab:sim}, the argument can be put to
work in the case $\alpha = 1/2$ that corresponds to median and mean
forecasts, respectively.  Specifically, let $F$ be the perfect
forecast, which has median and mean $\mu$, let $\cA$ be the sigma
field generated by $\mu$, and let $X_1 = 0$ and $X_2 = - \mu$.
Invoking the order sensitivity argument, we see that the
climatological forecaster dominates the sign-reversed forecaster for
both median and mean predictions.

\subsection{The Murphy diagram as a diagnostic tool}  \label{subsec:Murphy}

As noted, Corollaries 1a and 1b suggest graphical tools for the
comparison of quantile and expectile forecasts, including the special
cases of the mean or expectation functional, and the further special
case of probability forecasts of a binary event.  We describe these
diagnostic tools in the setting of a point prediction space
(\ref{eq:Xspace}), where $X_1, \ldots, X_l$ denote point forecasts for
the outcome $Y$, and the probability measure $\myQ$ represents their
joint distribution.  In the case of probability forecasts, we use the
more suggestive notation $p_1, \ldots, p_l$ for the forecasts.

\begin{itemize} 
\item 
For quantile forecasts at level $\alpha \in (0,1)$, we plot the graph
of the expected elementary quantile score $\mySQ$,
\begin{equation}  \label{eq:sQ}
\theta \mapsto s_j(\theta) = \myE_\myQ \hsp \mySQ(X_j,Y), 
\end{equation} 
for $j = 1, \ldots, l$.  By Corollary 1a, forecast $X_i$ dominates
forecast $X_j$ if and only if $s_i(\theta) \leq s_j(\theta)$ for
$\theta \in \real$.  The area under $s_j(\theta)$ equals the
respective expected asymmetric piecewise linear score (\ref{eq:apl}).
\item 
For expectile forecasts at level $\alpha \in (0,1)$, we plot the graph
of the expected elementary expectile score $\mySE$,
\begin{equation}  \label{eq:sE}
\theta \mapsto s_j(\theta) = \myE_\myQ \hsp \mySE(X_j,Y), 
\end{equation} 
for $i = 1, \ldots, l$.  By Corollary 1b, forecast $X_i$ dominates
forecast $X_j$ if and only if $s_i(\theta) \leq s_j(\theta)$ for
$\theta \in \real$.  The area under $s_j(\theta)$ equals half the
respective expected asymmetric squared error (\ref{eq:ase}).
\item
For probability forecasts of a binary event, we plot the graph
of the expected elementary score $\mySB$,
\begin{equation}  \label{eq:sB}
\theta \mapsto s_j(\theta) = \myE_\myQ \hsp \mySB( \, p_j,Y), 
\end{equation} 
for $i = 1, \ldots, l$.  By Corollary 1b, the probability forecast
$p_i$ dominates $p_j$ if and only if $s_i(\theta) \leq s_j(\theta)$
for $\theta \in (0,1)$.  The area under $s_j(\theta)$ equals half the
expected Brier score (\ref{eq:Brier}).  
\end{itemize}

In the context of probability forecasts for binary weather events,
displays of this type have a rich tradition that can be traced to
Thompson and Brier (1955) and Murphy (1977).  More recent examples
include the papers by Schervish (1989), Richardson (2000), Wilks
(2001), Mylne (2002), and Berrocal et al.~(2010), among many others.
Murphy (1977) distinguished three kinds of diagrams that reflect the
economic decisions involved.  The negatively oriented {\em expense
diagram}\/ shows the mean raw loss or expense of a given forecast
scheme; the positively oriented {\em value diagram}\/ takes the
unconditional or climatological forecast as reference and plots the
difference in expense between this reference forecast and the forecast
at hand, and lastly, the {\em relative-value diagram}\/ plots the
ratio of the utility of a given forecast and the utility of an oracle
forecast.  The displays introduced above are similar to the value
diagrams of Murphy, and we refer to them as {\em Murphy diagrams}.
Our Murphy diagrams are by default negatively oriented and plot the
expected elementary score for competing quantile, expectile, and
probability forecasters.  For better visual appearance, we generally
connect the left- and right-hand limits at the jump points of the
empirical score curves.

\begin{figure}[t]

\begin{minipage}{0.5\columnwidth}
\hspace*{33.5mm} \small Mean \\
\includegraphics[width=\columnwidth]{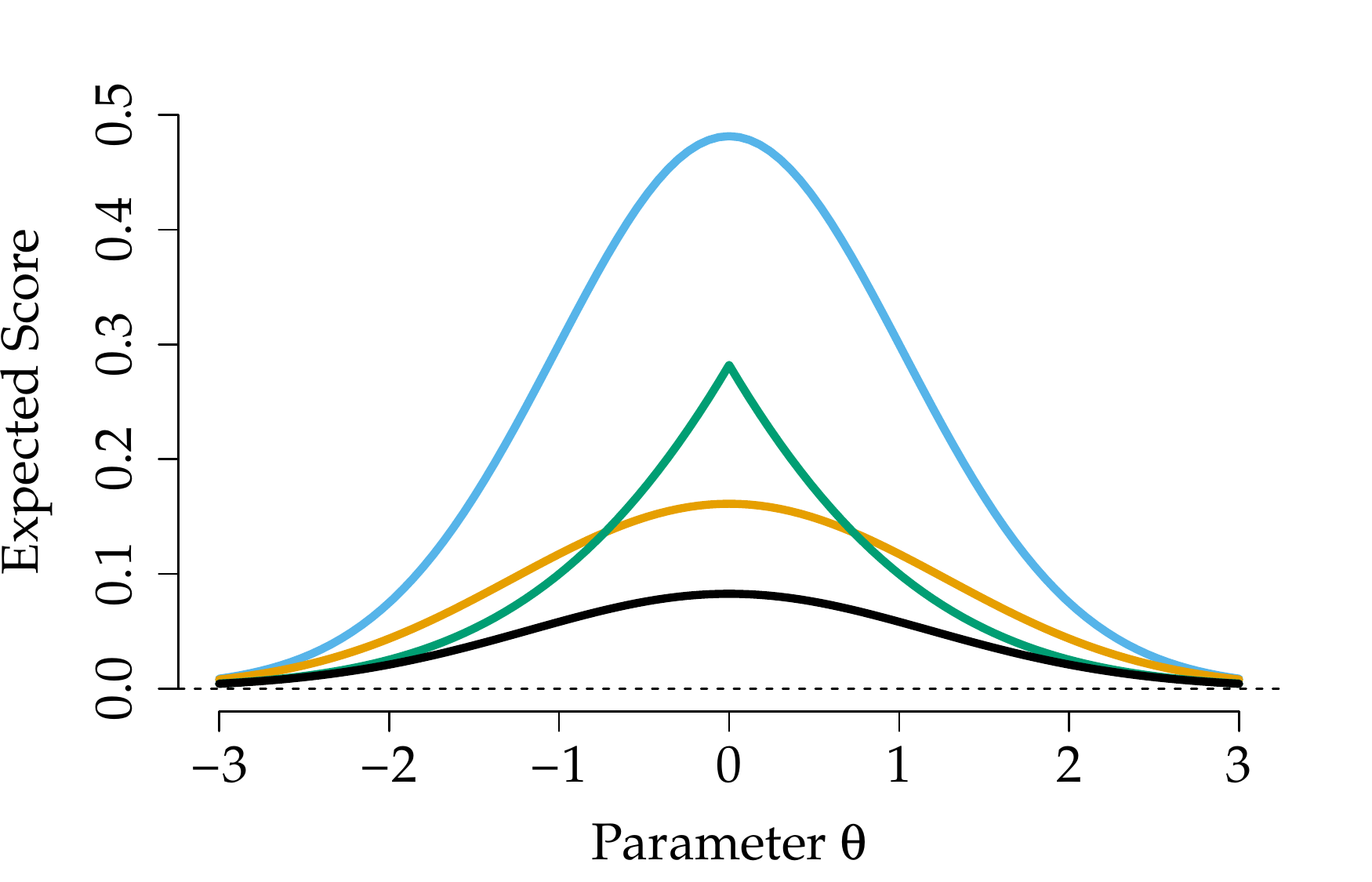}
\hspace*{20mm} \small Quantile ($\alpha = 0.90$) \rule{0mm}{7mm} \\
\includegraphics[width=\columnwidth]{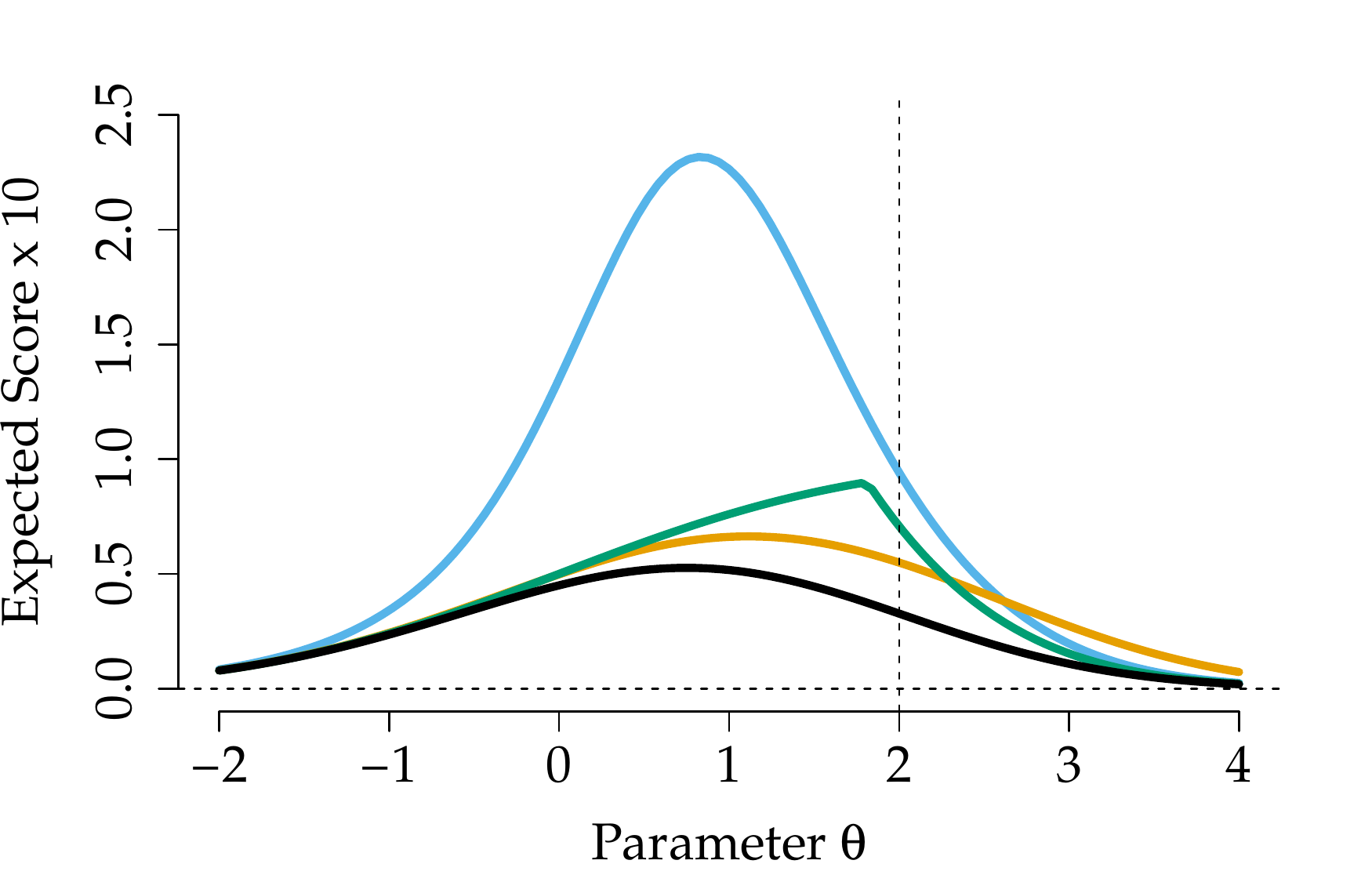}
\end{minipage}
\begin{minipage}{0.5\columnwidth}
\hspace*{25mm} \small Mean \\
\includegraphics[width=\columnwidth]{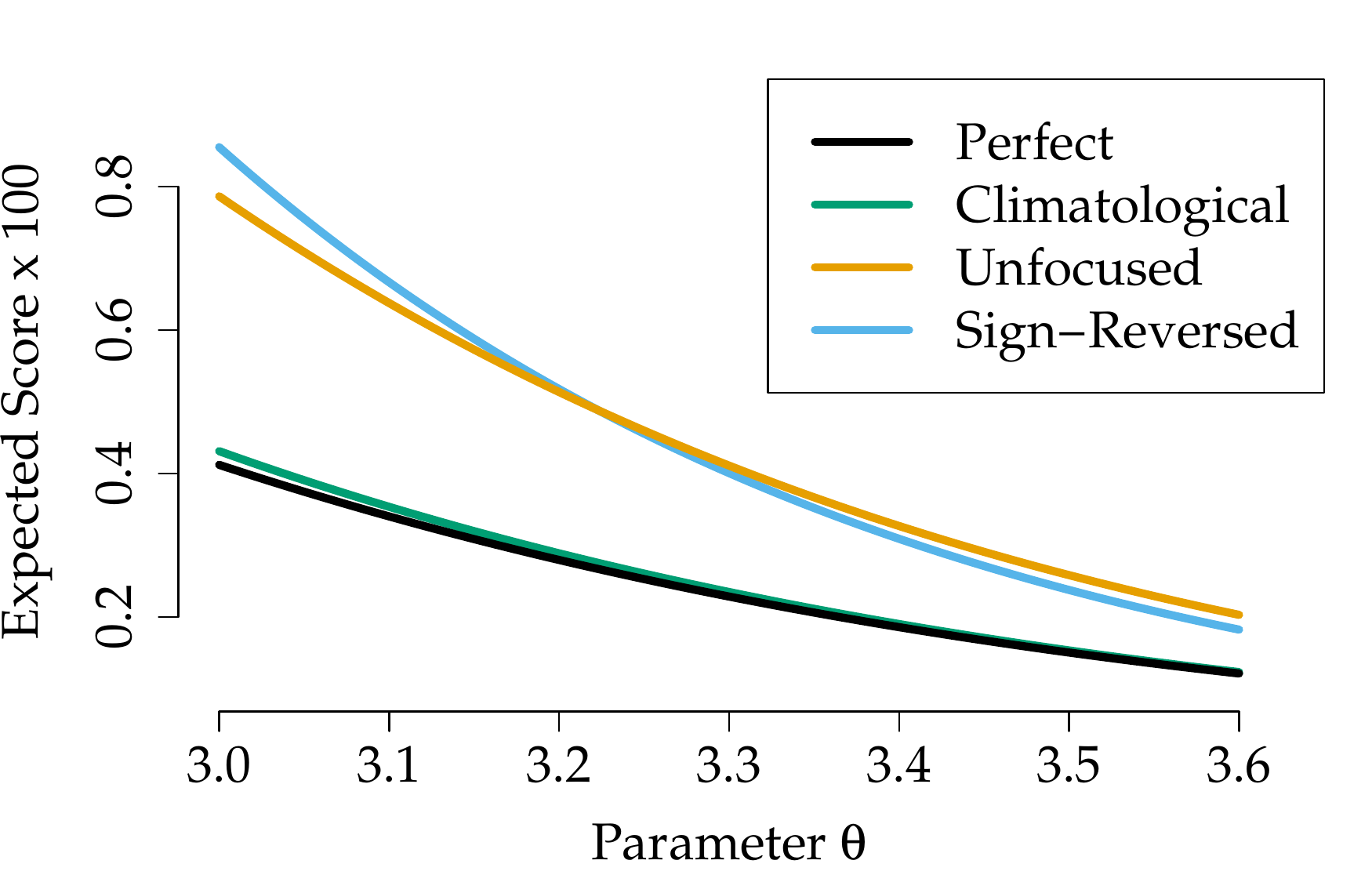}
\hspace*{20mm} \small Probability ($Y \geq 2$) \rule{0mm}{7mm} \\
\includegraphics[width=\columnwidth]{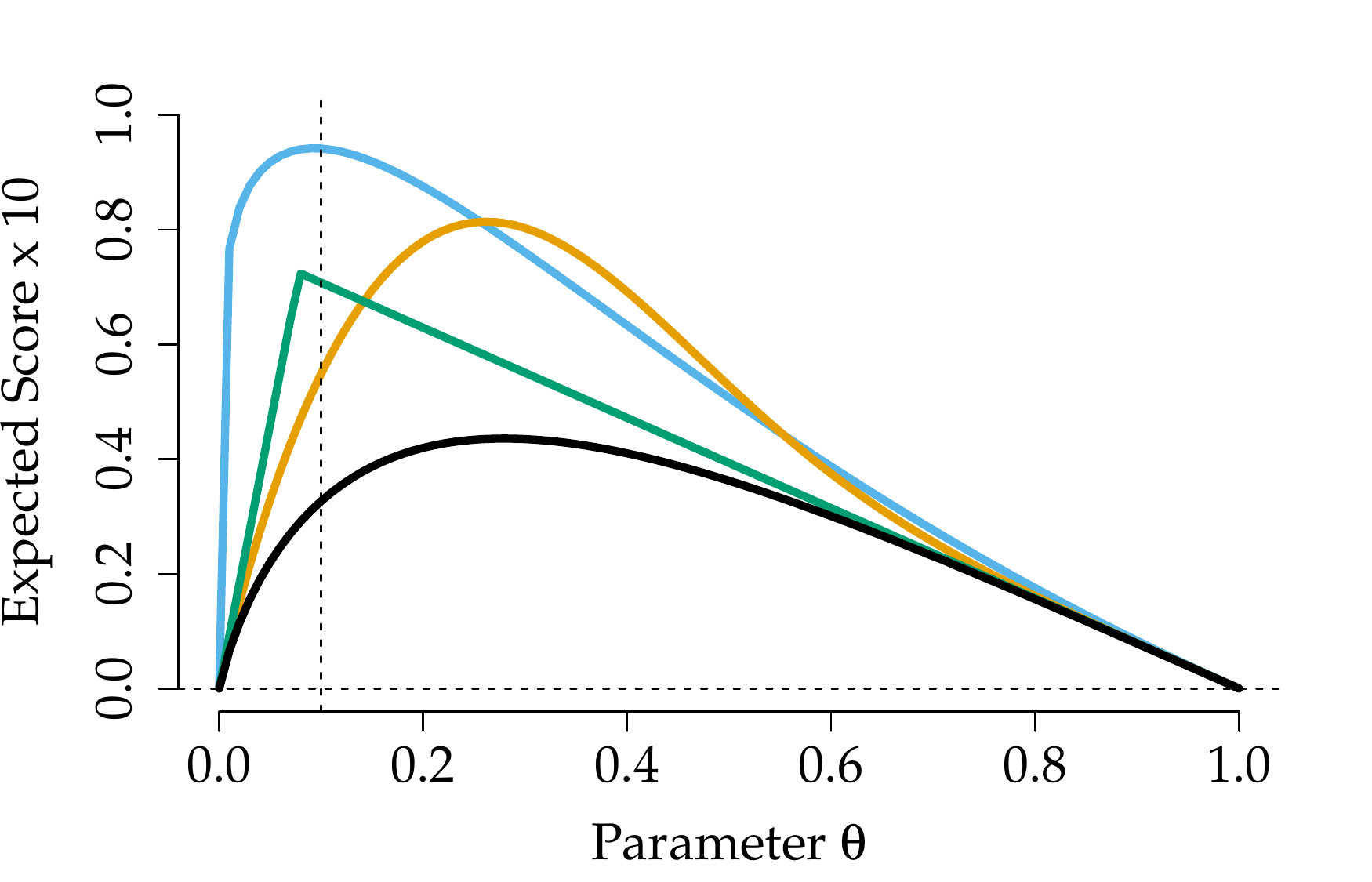}
\end{minipage}

\caption{\small Murphy diagrams for the forecasters in Table
  \ref{tab:sim}.  The functionals considered are the mean, the
  quantile at level $\alpha = 0.90$, and the probability of the binary
  event $Y \geq 2$.  The vertical dashed lines in the bottom panels
  indicate the extremal scores $\mySQ$ and $\mySB$ that relate to each
  other as in (\ref{eq:psfx}) and (\ref{eq:qep}).
  \label{fig:sim}}

\end{figure}

Figure \ref{fig:sim} shows Murphy diagrams for the perfect,
climatological, unfocused, and sign-reversed forecasters in Table
\ref{tab:sim}.  We compare point predictions for the mean or
expectation functional, and the quantile at level $\alpha = 0.90$,
along with probability forecasts for the binary event that the outcome
exceeds the threshold value 2.  Analytic expressions for the
respective expected scores are given in Appendix B.  As proved in the
previous section, the perfect forecaster dominates the other
forecasters for all functionals considered.  The expected score curves
for the climatological and the unfocused, and for the unfocused and
the sign-reversed forecasters, intersect in all three cases, so there
are no order relations between these forecasters.  Finally, the Murphy
diagrams suggest that the climatological forecaster dominates the
sign-reversed forecaster for all three functionals, and in the case of
the mean functional, the order sensitivity argument in the previous
section confirms the visual impression.  In the cases of the quantile
and probability forecasts, final confirmation would need to be based
on tedious analytic investigations of the asymptotic behavior of the
expected score functions.

By default, our Murphy diagrams show the expected elementary scores.
If interest focuses on binary comparisons, it is natural to
consider Murphy diagrams for the difference,
\begin{equation}  \label{eq:D} 
\theta \mapsto D(\theta) = 
\myE_\myQ \hsp \myS_\theta(X_1,Y) - \myE_\myQ \hsp \myS_\theta(X_2,Y).
\end{equation} 
between the expected elementary scores of two point forecasters.

\subsection{Murphy diagrams for empirical forecasters}  \label{subsec:empirical} 

We now turn to the comparison and ranking of empirical forecasts.
Specifically, we consider tuples
\begin{equation}  \label{eq:empirical} 
\left( x_{i1}, \ldots, x_{il}, y_i \right),  \qquad i = 1, \ldots, n, 
\end{equation} 
where $x_{1j}, \ldots, x_{nj}$ are the $j$th forecaster's point
predictions, for $j = 1, \ldots, l$, and $y_i, \ldots, y_n$, are the
respective outcomes.  Thus, we have $l$ competing forecasters, and
each of them issues a set of $n$ point predictions.  A convenient
interpretation of the empirical setting is as a special case of a
point prediction space, in which the tuples $(X_1, \ldots, X_l, Y)$ in
(\ref{eq:Xspace}) attain each of the values in (\ref{eq:empirical})
with probability $1/n$.  Then the probability measure $\myQ$ is the
corresponding empirical measure, and with this identification, 
the (average) empirical scores
\[
s_j(\theta) = \frac{1}{n} \sum_{i=1}^n \myS_\theta(x_{ij},y_i), 
\]
where $\myS_\theta$ is either $\mySQ$, $\mySE$, or $\mySB$, become the
expected elementary scores from (\ref{eq:sQ}), (\ref{eq:sE}), and
(\ref{eq:sB}), respectively.  To compare forecasters $X_1$ and $X_2$,
say, it is convenient to show a Murphy plot of the equivalent of the
difference (\ref{eq:D}), namely
\[
\theta \mapsto D_n(\theta) = \frac{1}{n} \sum_{i=1}^n \, d_i(\theta), 
\]
where
\begin{equation}  \label{eq:d}
d_i(\theta) = \myS_\theta(x_{i1},y_i) - \myS_\theta(x_{i2},y_i)
\end{equation} 
for $i = 1, \ldots, n$, and again $\myS_\theta$ is either $\mySQ$,
$\mySE$, or $\mySB$, respectively.

Murphy diagrams can be used efficiently to show a lack of domination
when forecasters' expected elementary score curves intersect.
However, in general it is not possible to conclude domination, unless
the visual impression is supported by tedious analytic investigations
of the behavior of the expected score functions as $\theta \to \pm
\infty$.  Fortunately, these complications do not arise in the
empirical case, where dominance can be established by comparing the
empirical score functions at a well-defined, finite set of
arguments only, as follows.

\medskip \noindent 
{\bf Corollary 2a (quantiles).}  
{\em An empirical forecast\/ $X_1$ dominates\/ $X_2$ for
  $\alpha$-quantile predictions if\/ 
\[
\frac{1}{n} \sum_{i=1}^n \mySQ(x_{1i},y_i) \leq \frac{1}{n} \sum_{i=1}^n \mySQ(x_{2i},y_i)
\] 
for\/ $\theta \in \{ x_{11}, x_{12}, y_1, \ldots, x_{n1}, x_{n2}, y_n \}$.}

\medskip \noindent
{\bf Corollary 2b (expectiles).}  
{\em An empirical forecast\/ $X_1$ dominates\/ $X_2$ for
  $\alpha$-expectile predictions if\/ 
\[
\frac{1}{n} \sum_{i=1}^n \mySE(x_{1i},y_i) \leq 
\frac{1}{n} \sum_{i=1}^n \mySE(x_{2i},y_i)
\] 
for\/ $\theta \in \{ x_{11}, x_{12}, y_1, \ldots, x_{n1}, x_{n2}, y_n
\}$ and in the left-hand limit as\/ $\theta \uparrow \theta_0 \in \{
x_{11}, x_{12}$, $\ldots, x_{n1}, x_{n2} \}$.  In the case\/ $\alpha =
1/2$ evaluations at\/ $\theta \in \{ y_1, \ldots, y_n \}$ can be
omitted.}

\begin{figure}[t]

\centerline{
\includegraphics[width=0.50\columnwidth]{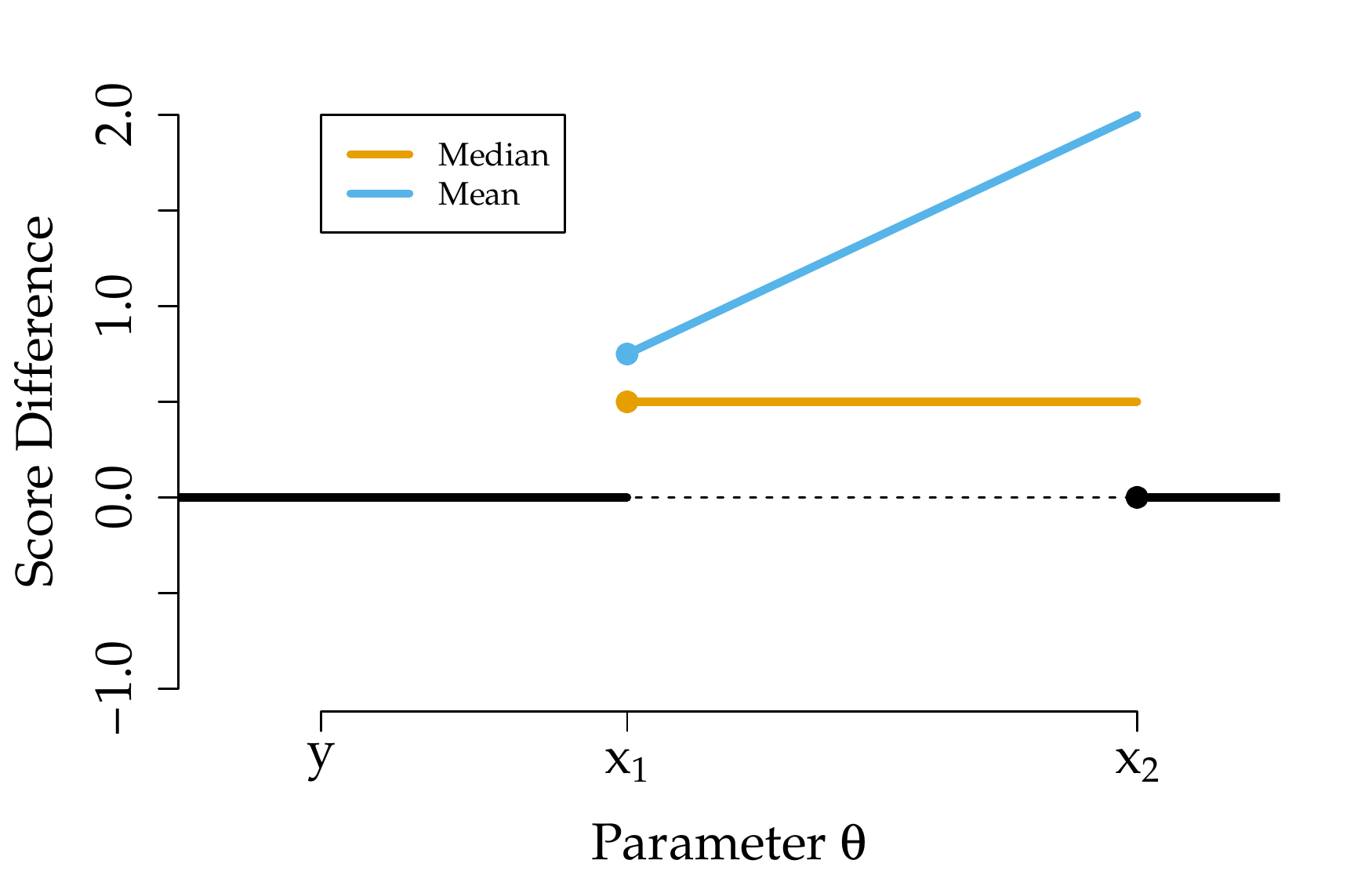}
\includegraphics[width=0.50\columnwidth]{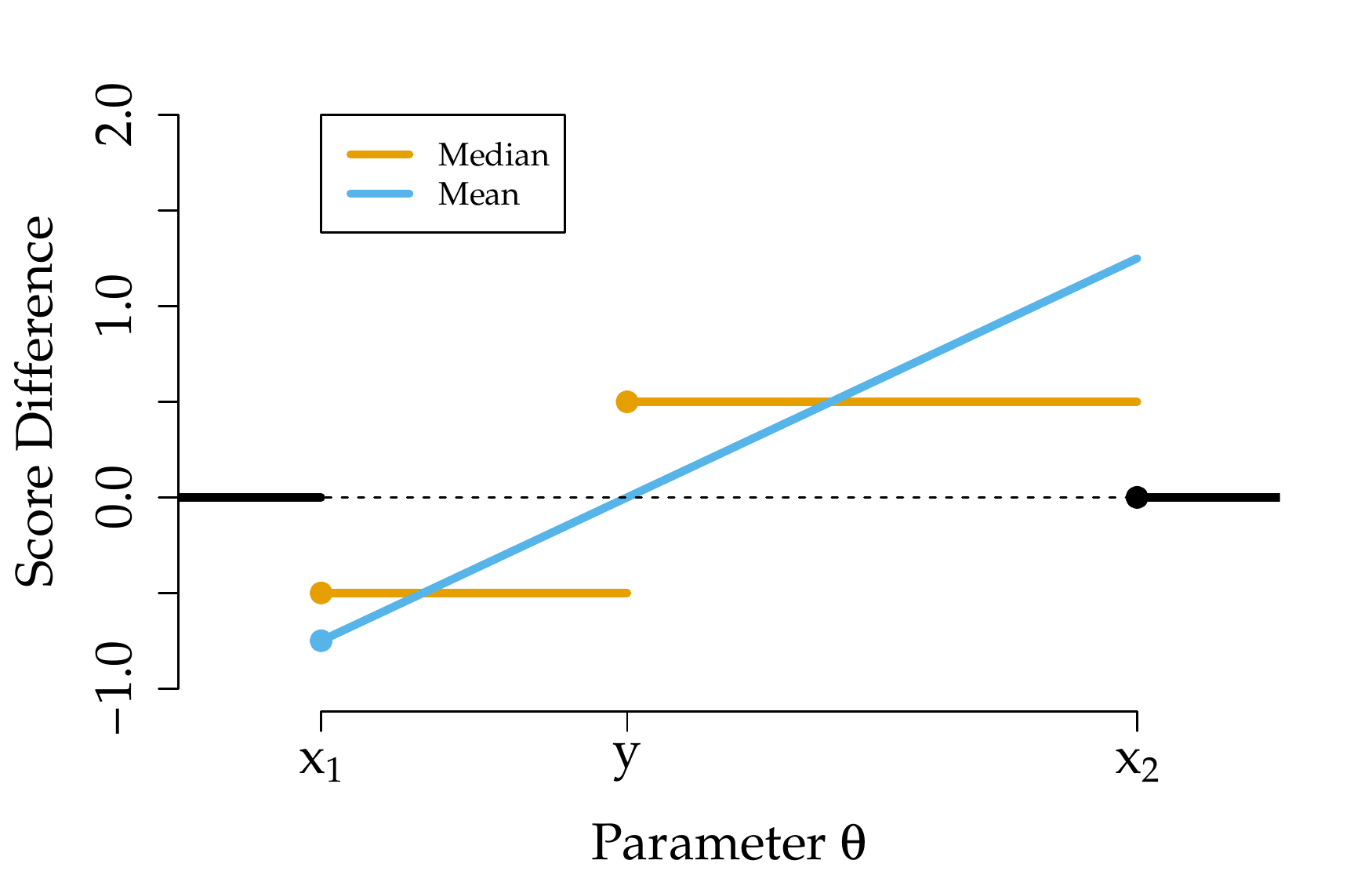}
}

\caption{\small The general shape of the score differential
  $d_i(\theta)$ in (\ref{eq:d}) for the median and mean
  functionals. \label{fig:steps}}

\end{figure}

\medskip 
To see why these results hold, note that in either case the score
differential $d_i(\theta)$ is right-continuous, and that it vanishes
unless $ \min(x_{i1},x_{i2}) \leq \theta < \max(x_{i1},x_{i2})$.
Furthermore, in the case of quantiles $d_i(\theta)$ is piecewise
constant with no other jump points than $x_{i1}, x_{i2}$, or $y_i$.
Similarly, in the case of expectiles $d_i(\theta)$ is piecewise linear
with no other jump points than $x_{i1}$ and $x_{i2}$, and no other 
change of slope than at $y_i$.  The change of slope disappears when
$\alpha = 1/2$.  Figure \ref{fig:steps} illustrates the behavior of
$d_i(\theta)$ in the cases of the median and the mean, respectively.

\begin{figure}[t]

\centerline{
\includegraphics[width=0.50\columnwidth]{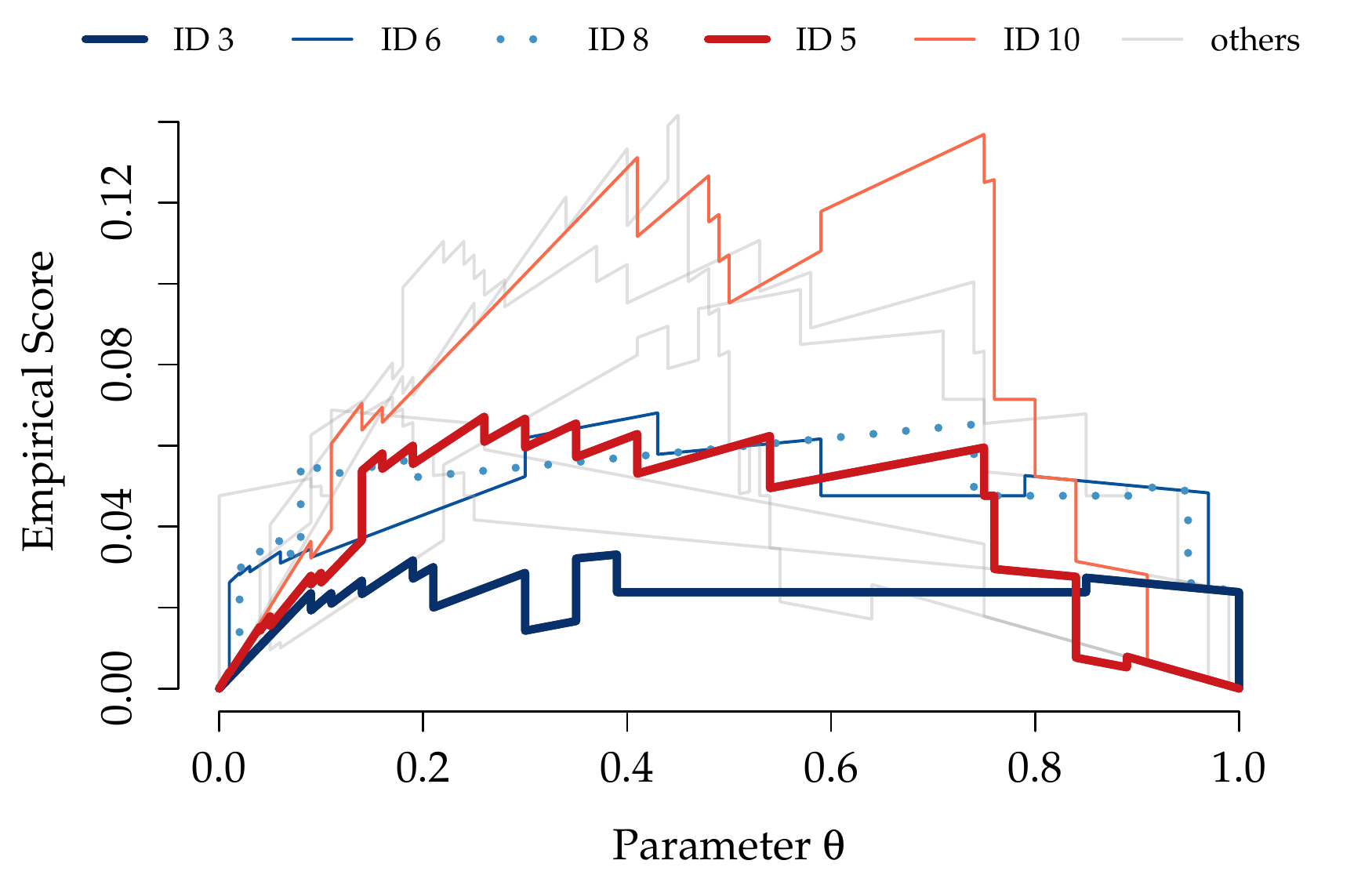}
\includegraphics[width=0.50\columnwidth]{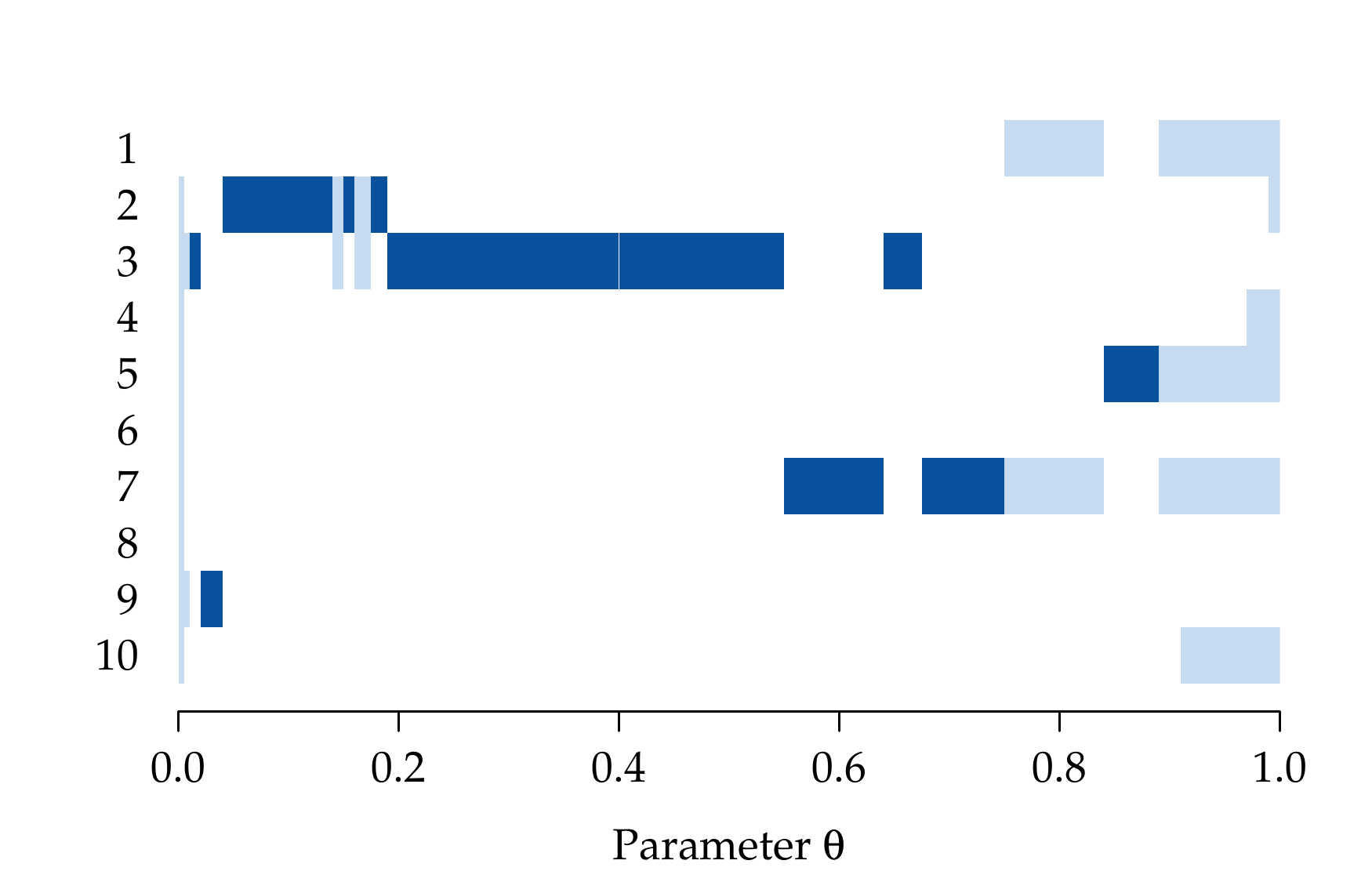}
}

\caption{\small Left: Murphy diagram for the probability forecasters
  in Table A.1 of Merkle and Steyvers (2013).  Right: The best
  forecast ID(s) under $\mySB$, with dark blue indicating a unique
  best score, and light blue a shared best score.  For example, ID 9
  attains the unique best score for $\theta \in [0.02, 0.04)$, and ID
    10 attains the shared best score for $\theta \in [0.91,
      1)$.  \label{fig:ACES}}

\end{figure}

To give an example, we consider the 10 forecasters in Table A.1 of
Merkle and Steyvers (2013), each of whom issues probability forecasts
for 21 binary events.  The data are artificial but mimic forecasters
in the Aggregate Contingent Estimation System (ACES), a web based
survey that solicited probability forecasts for world events from the
general public.  The Murphy diagram in the left-hand panel of Figure
\ref{fig:ACES} shows the empirical score curves
\[
\theta \mapsto s_j(\theta) = \frac{1}{21} \sum_{i=1}^{21} \mySB( \hsp p_{ij},y_i), 
\]
where $p_{ij} \in [0,1]$ is forecaster $j$'s stated probability for
world event $i$ to materialize, and $y_i \in \{ 0, 1 \}$ is the
respective binary realization.  By Corollary 2b, dominance relations
can be inferred by evaluating $s_j(\theta)$ at the forecasters' stated
probabilities.  We note that ID 3 dominates IDs 6 and 8, and that ID 5
dominates ID 10.  The remaining pairwise comparisons do not give rise
to dominance relations.  The induced partial order between the IDs
applies to comparisons under {\em any}\/ proper scoring rule, as
reflected by the rankings in Table 1 of Merkle and Steyvers (2013).
The right-hand panel in Figure \ref{fig:ACES} considers joint
comparisons.  We see that ID 3 attains the lowest score over a wide
range of $\theta$.  However, IDs 2, 5, 7, and 9 show the unique best
empirical score under $\mySB$ for other values of $\theta$ and,
therefore, have superior economic utility under the associated
cost-loss ratios.

\section{Empirical examples}  \label{sec:examples}

We now demonstrate the use of Murphy diagrams in economic and
meteorological case studies in time series settings.  In each example,
interest is in a comparison of two forecasts, and so we show Murphy
diagrams for the empirical scores and their difference.  The jagged
visual appearance stems from the behavior of the empirical score
functions just explained and depends on the number $n$ of forecast
cases.  We supplement the Murphy diagrams for a difference by
confidence bands based on Diebold and Mariano (1995) tests with a
heteroscedasticity and autocorrelation robust variance estimator
(Newey and West 1987).  The approach of Diebold and Mariano (1995)
views empirical data of the form (\ref{eq:empirical}) as a sample from
an underlying population and tests the hypothesis of equal expected
scores.  The confidence bands are pointwise and have a nominal level
of 95\%.

\begin{figure}[p]

\centering

\begin{tabular}{c} 
\toprule
\text{\small Mean Inflation, Patton (2015); $n = 129$ \rule{0mm}{7mm}} \vspace{-7mm} \\ 
\includegraphics[width = \textwidth]{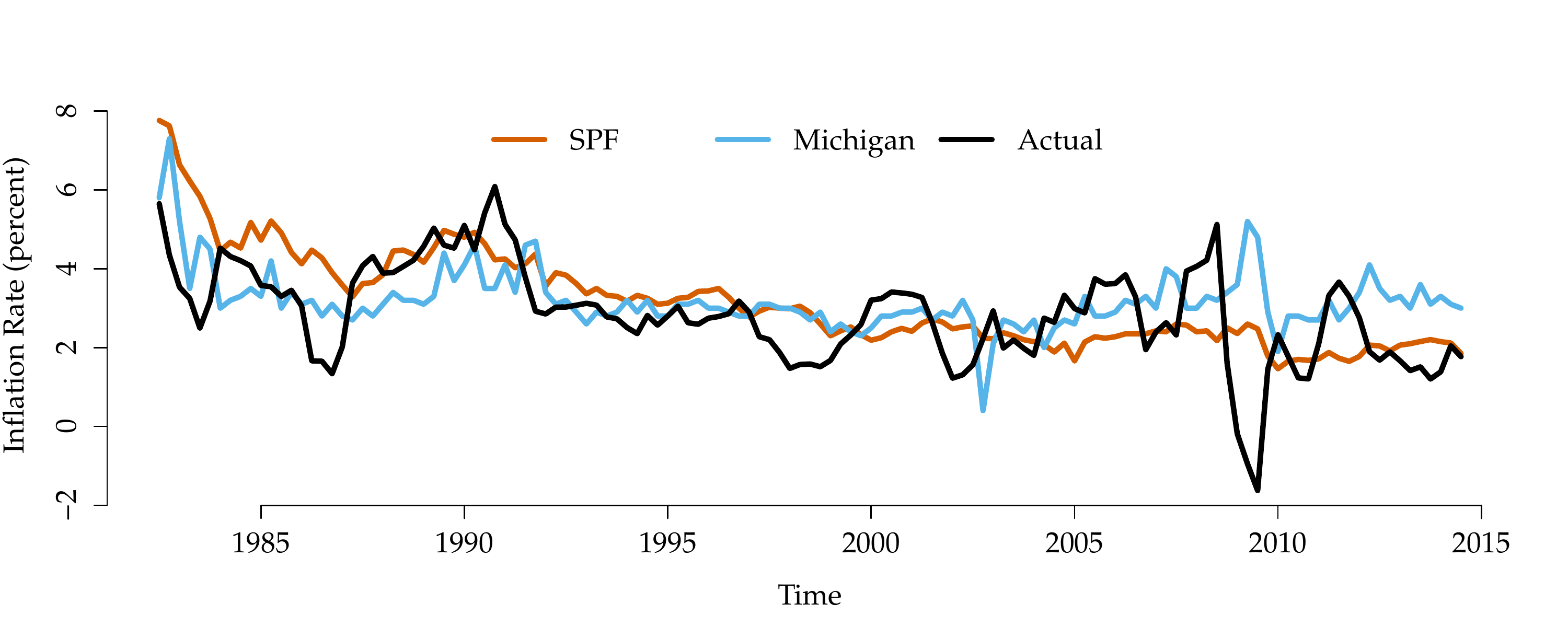} \\ 
\text{\small Probability of Recession, Rudebusch and Williams (2009); $n = 186$ \rule{0mm}{7mm}} \vspace{2mm} \\ 
\includegraphics[width = \textwidth]{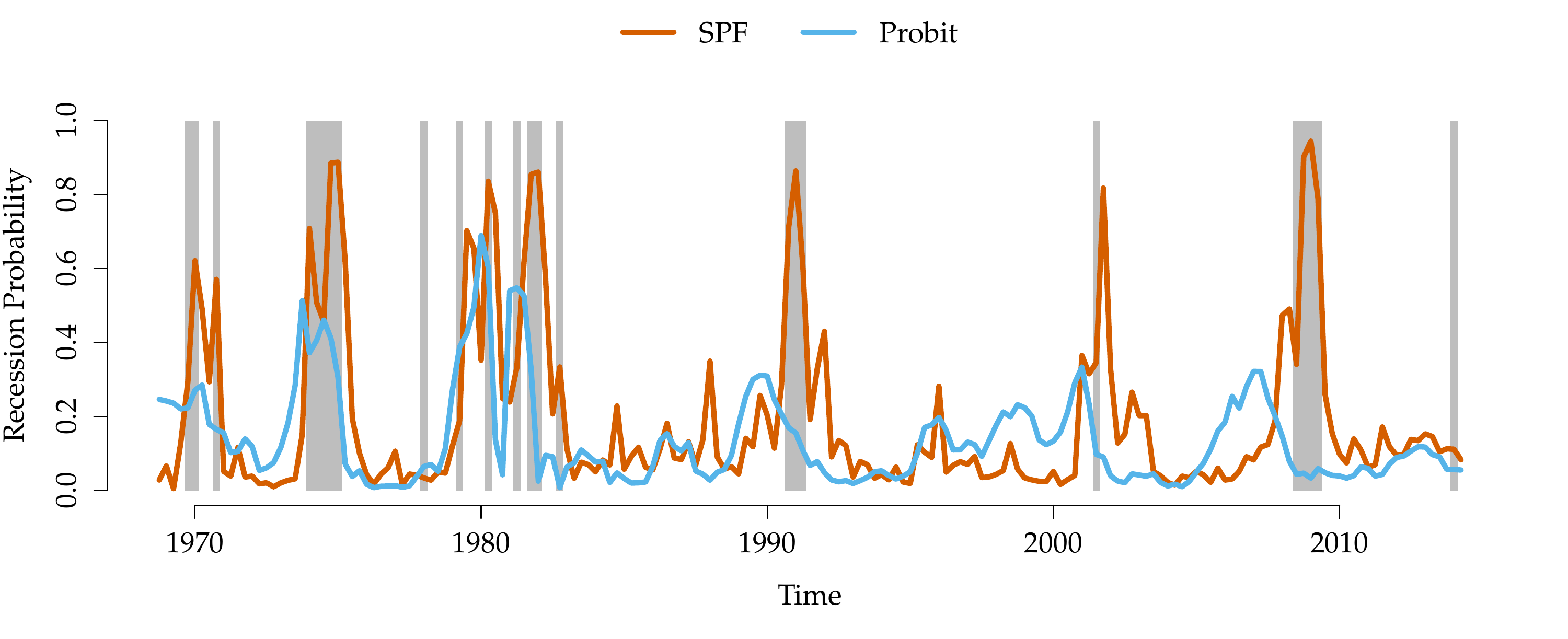} \\
\text{\small 90\% Quantile of Wind Speed, Gneiting et al.~(2006); $n = 5136$ \rule{0mm}{7mm}} \vspace{-7mm} \\ 
\includegraphics[width = \textwidth]{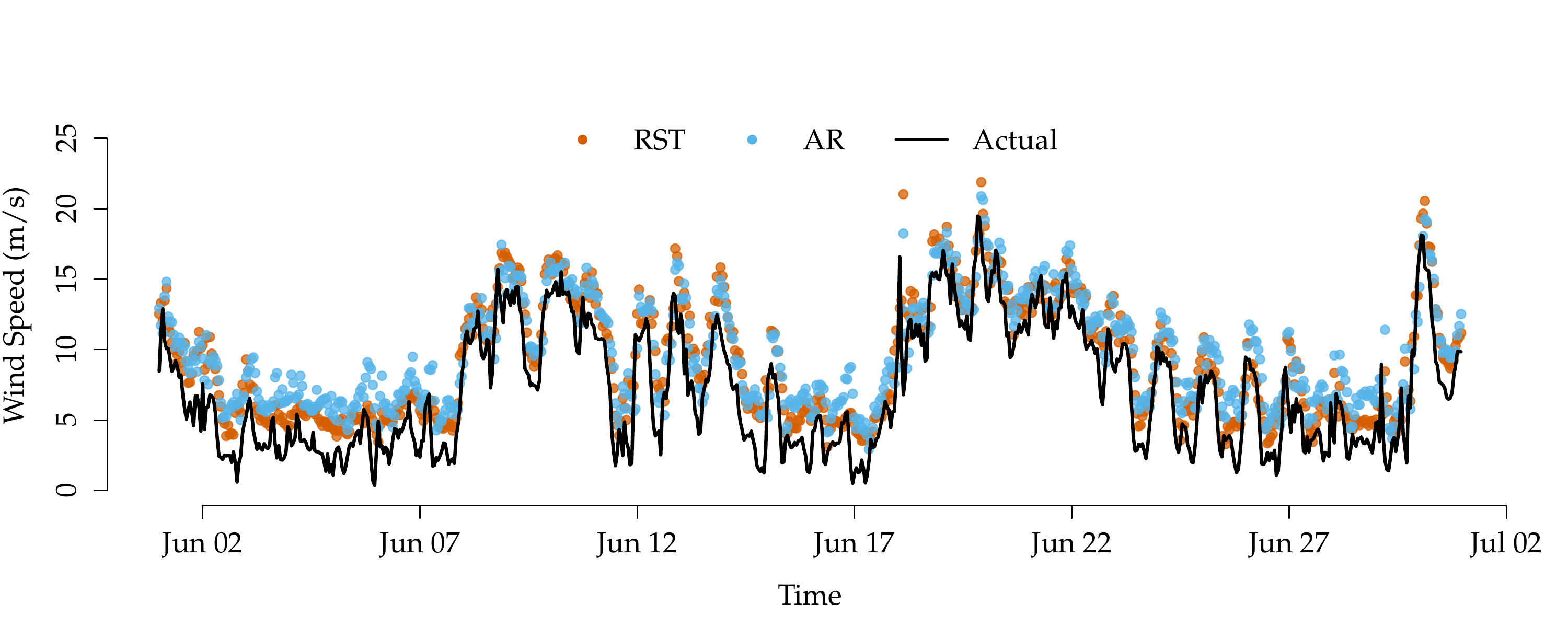} \\ 
\bottomrule
\end{tabular}

\caption{\small Point forecasts and realizations in the empirical
  examples.  In the middle plot, shaded areas correspond to actual
  recessions.  The plot at bottom is restricted to a subperiod in summer 
  2003.  \label{fig:data}}

\end{figure}

\begin{figure}[p]
	
\centering

\[
\begin{array}{cc}
\toprule
\text{\small Score} & \text{\small Score Difference \rule{0mm}{4mm}} \\
\midrule
\multicolumn{2}{c}{\text{\small Mean Inflation, Patton (2015); $n = 129$ \rule{0mm}{7mm}}} \\ 
\includegraphics[height = 4.5cm]{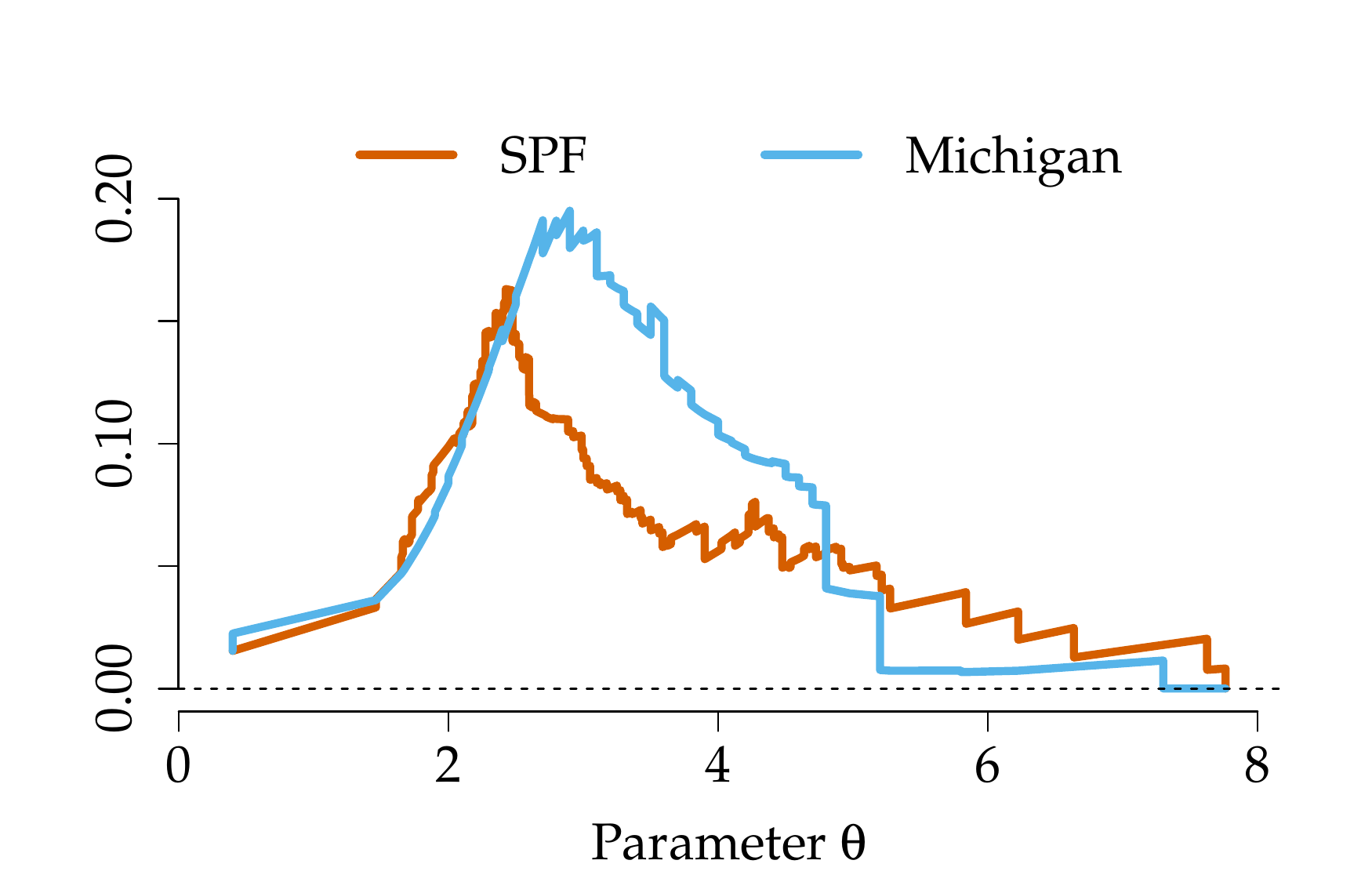} &
\includegraphics[height = 4.5cm]{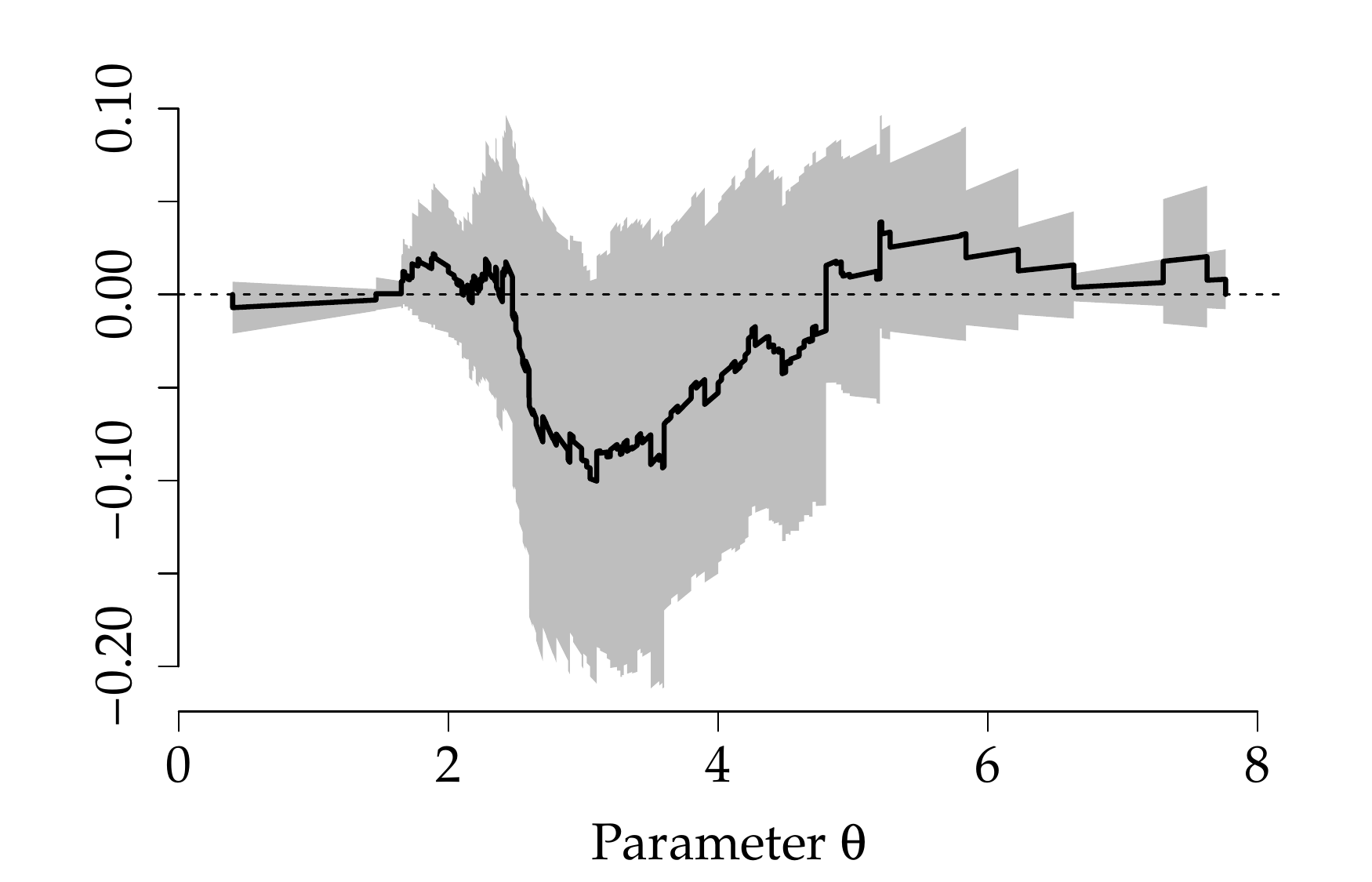} \\ 
\multicolumn{2}{c}{\text{\small Probability of Recession, Rudebusch and Williams (2009); $n = 186$ \rule{0mm}{7mm}}} \\ 
\includegraphics[height = 4.5cm]{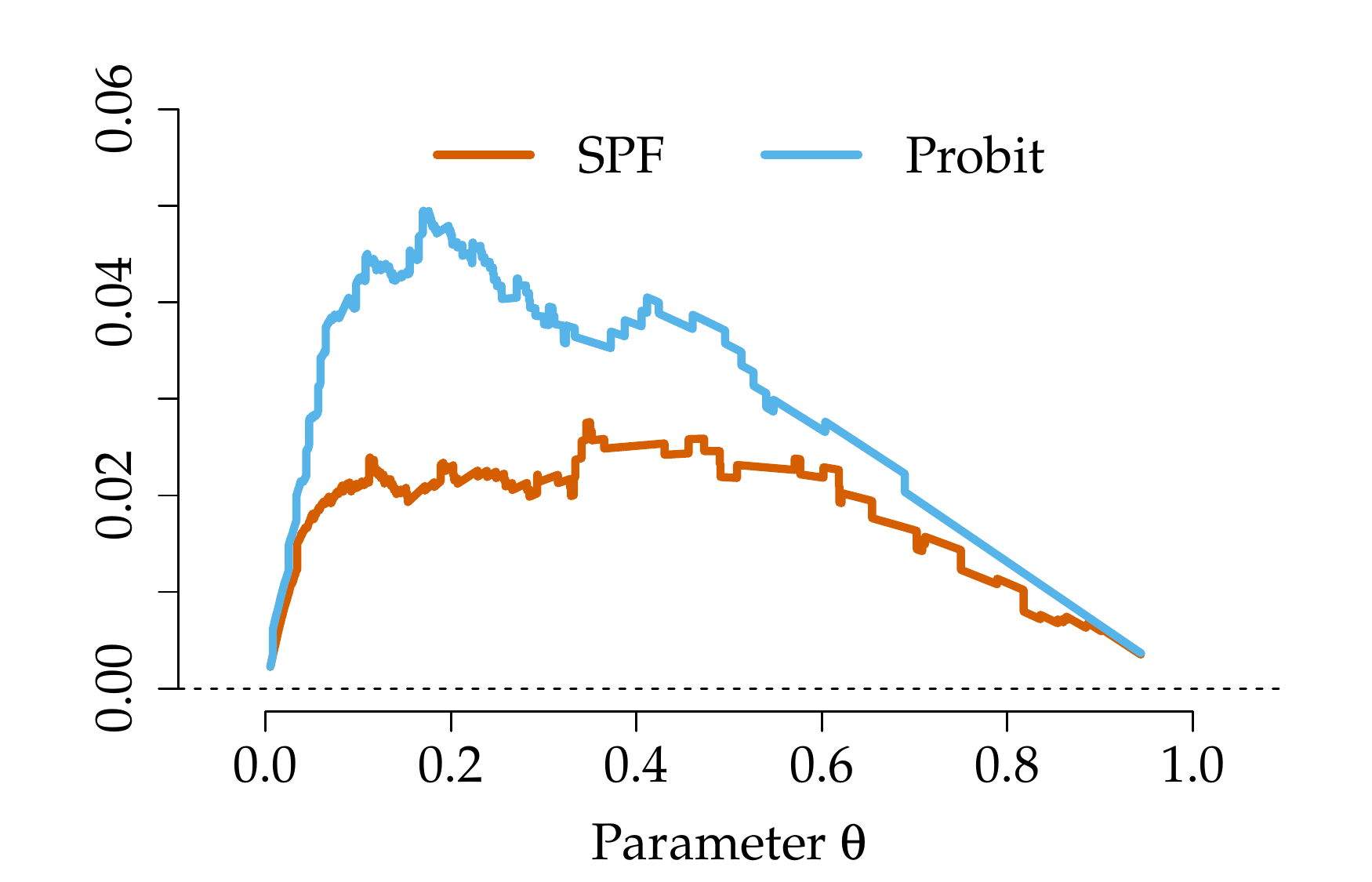} &
\includegraphics[height = 4.5cm]{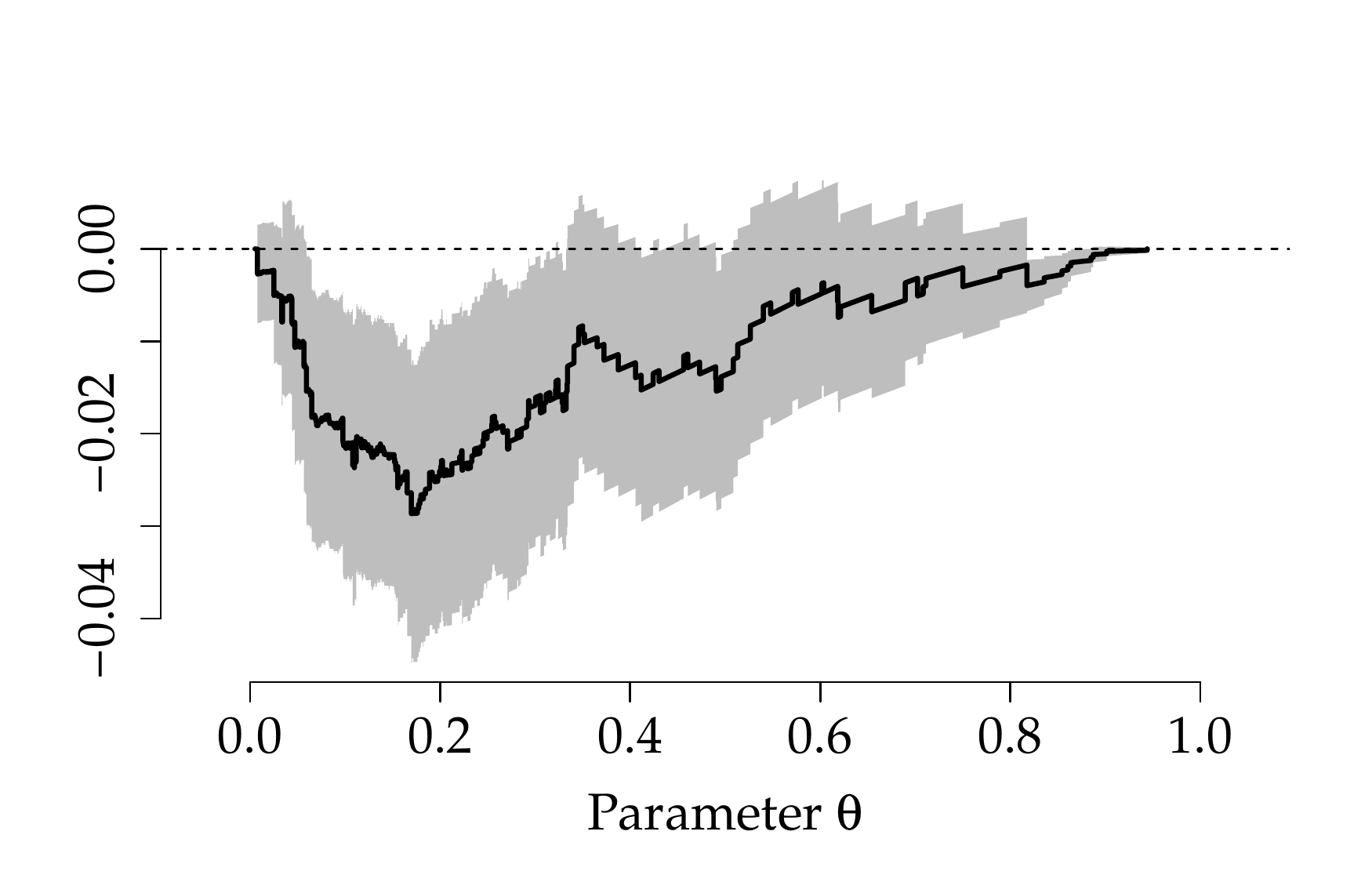} \\ 
\multicolumn{2}{c}{\text{\small 90\% Quantile of Wind Speed, Gneiting et al.~(2006); $n = 5136$ \rule{0mm}{7mm}}} \\ 
\includegraphics[height = 4.5cm]{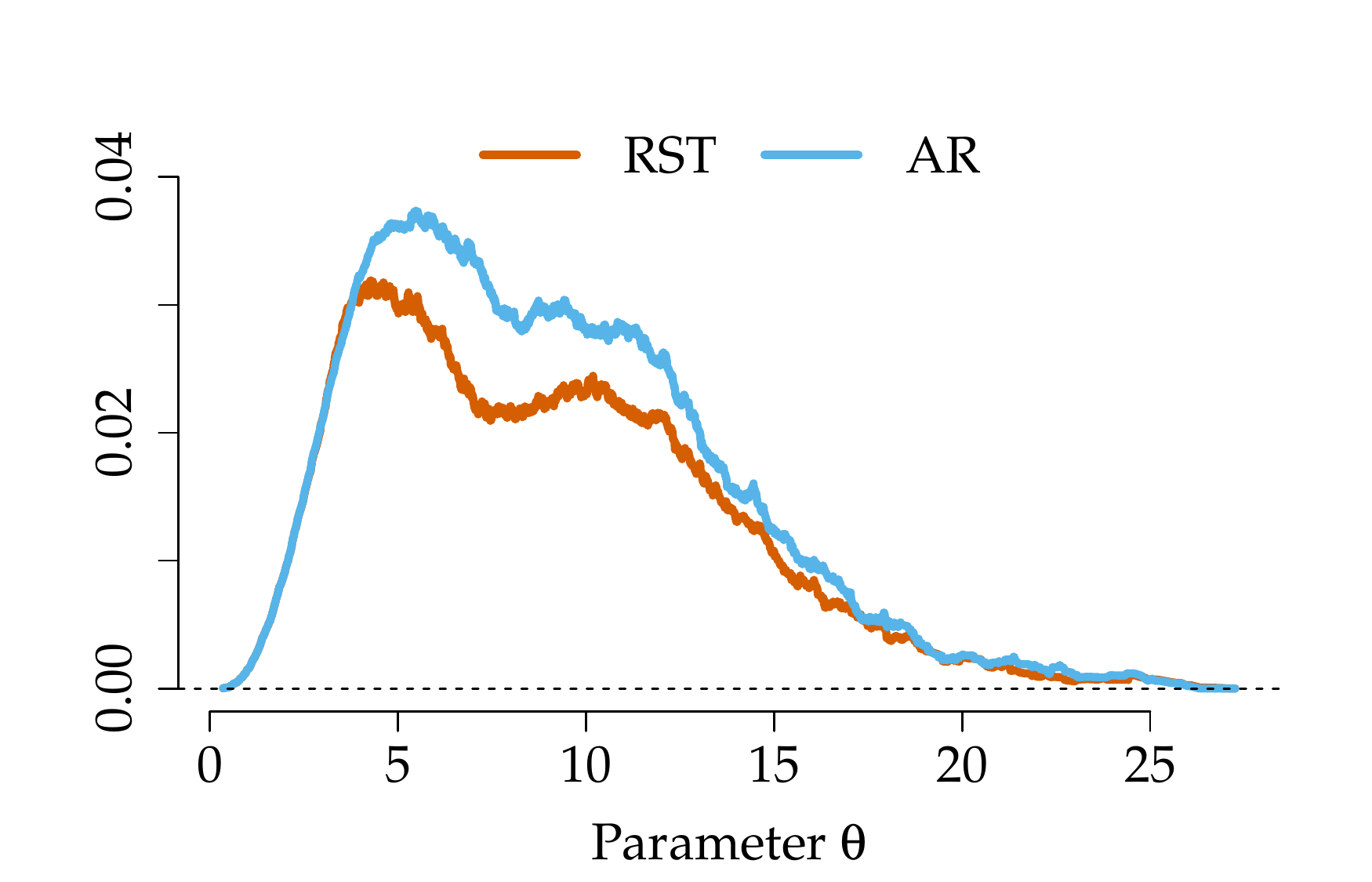} &
\includegraphics[height = 4.5cm]{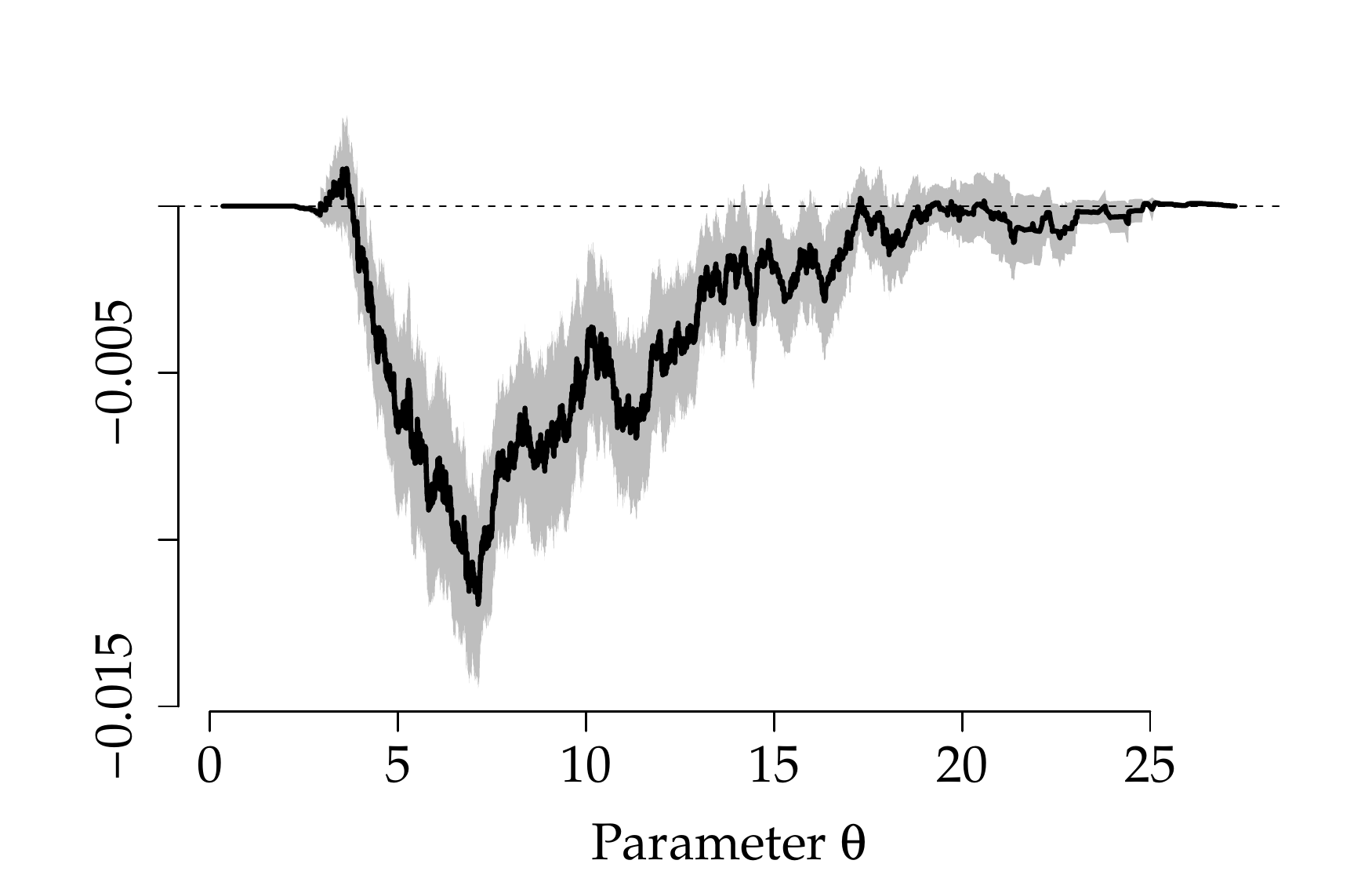} \\ 
\bottomrule
\end{array}
\]

\caption{\small Murphy diagrams in the empirical examples.  In the right
  column, a negative score difference means that the first named
  method is preferable.  \label{fig:Murphy}} 

\end{figure}

\subsection{Mean forecasts of inflation}  \label{subsec:SPF} 

In macroeconomics, subjective expert forecasts often compare favorably
to statistical forecasting approaches; see Faust and Wright (2013) for
evidence and discussion.  For the United States, the Survey of
Professional Forecasters (SPF) run by the Federal Reserve Bank of
Philadelphia is a key data source; see, e.g., Engelberg { et
  al.}~(2009). Patton (2015) uses SPF data to illustrate the use of
various scoring functions that are consistent for the mean functional.

Motivated by Patton's analysis, we analyze SPF mean forecasts for the
annual inflation rate of the Consumer Price Index (CPI).  We compare
the SPF forecasts to forecasts from another survey, the Michigan
Survey of Consumers, based on data from the third quarter of 1982 to
the third quarter of 2014, for a test period of 129 quarters.  Our
implementation choices are as in Section 5 of Patton (2015), except
that we update the data set to cover the observations for the second
and third quarters in 2014, and that we use the slightly newer fourth
quarter of 2014 vintage for the CPI realizations.  The top panel of
Figure \ref{fig:data} shows the forecasts along with the realizing
values.

The respective Murphy diagrams are shown in the top panel of Figure
\ref{fig:Murphy}.  At left, the curves for the empirical elementary
score $\myS^{\rm E}_{1/2,\theta}$ of the SPF and the Michigan survey
intersect prominently, suggesting that neither of the two surveys
dominates the other.  Specifically, the SPF is preferred for smaller
values, whereas the Michigan survey is preferred for larger values of
$\theta$.  This may be explained by a series of high inflation rates
up until 1992, which were better matched by the Michigan survey than
by the SPF.  At right, the confidence bands for the score differences
are fairly broad and include zero for all values of $\theta$.

\subsection{Probability forecasts of recession}  \label{subsec:recession} 

We now relate to the rich literature on binary regression and
prediction and analyze probability forecasts of United States
recessions, as proxied by negative real gross domestic product (GDP)
growth.  The SPF covers probability forecasts for this event since the
fourth quarter of 1968.  Following Rudebusch and Williams (2009), we
compare current quarter probability forecasts from the SPF to
forecasts from a probit model based on the term spread, i.e., the
difference between long and short term interest rates.  We follow
Rudebusch and Williams (2009) in all data and implementation choices,
except that we update their sample through the second quarter of 2014,
for a test period of 186 quarters.  Detailed economic and/or
statistical justification of these choices can be found in the
original paper.

The middle row of Figure \ref{fig:data} shows the SPF and probit model
based probability forecasts for a recession, with the gray vertical
bars indicating actual recessions.  During recessionary periods, the
SPF tends to assign higher forecast probabilities than the probit
model.  Also, the SPF tends to assign lower forecast probabilities
during non-recessionary periods.  The respective Murphy diagrams in
the middle row of Figure \ref{fig:Murphy} show that the SPF attains
lower empirical elementary scores $\mySB$ at all thresholds $\theta
\in (0,1)$.  The confidence bands for the score differences exclude
zero for small values of the cost-loss ratio $\theta$ and confirm the
superiority of the SPF over the probit model for current quarter
forecasts.  This can partly be attributed to the fact that SPF
panelists have access to timely within-quarter information that is not
available to the probit model.  As demonstrated by Rudebusch and
Williams (2009), the relative performance of the probit model improves
at longer forecast horizons, where within-quarter information plays a
lesser role.

\subsection{Quantile forecasts for wind speed}  \label{subsec:wind} 

We return to the meteorological example in Figure \ref{fig:wind}, but
instead of the mean or expectation functional we now consider quantile
forecasts at level $\alpha = 0.90$.  We compare the regime-switching
space-time (RST) approach introduced by Gneiting et al.~(2006) to a
simple autoregressive (AR) benchmark for two-hour ahead forecasts of
hourly average wind speed at the Stateline wind energy center in the
Pacific Northwest of the United States.  The original paper refers to
the specifications considered here as RST-D-CH and AR-D-CH,
respectively.  This terminology indicates that the methods account for
the diurnal cycle and conditional heteroscedasticity.  The data set,
evaluation period, estimation and forecast methods for this example
are identical to those in Gneiting et al.~(2006), and we refer to the
original paper for detailed descriptions.  Both methods yield
predictive distributions, from which we extract the quantile
forecasts.  The evaluation period ranges from May 1 through November
30, 2003, for a total of 5,136 hourly forecast cases.

The bottom panel in Figure \ref{fig:data} shows the quantile forecasts
and realizations.  The quantile forecasts exceed the outcomes at about
the nominal level, at 89.7\% for the RST forecast and 90.9\% for the
AR forecast, respectively, indicating good calibration.  However, the
RST forecasts are sharper, in that the average forecast value over the
evaluation period is 9.2 meters per second, as compared to 9.7 meters
per second in the case of the AR forecast.  To see why the sharpness
interpretation applies here, note that wind speed is a nonnegative
quantity, so quantile forecasts can be identified with one-sided
prediction intervals with a left limit of zero.  These observations
suggest the superiority of the RST forecasts over the benchmark AR
forecasts, and the Murphy diagrams for the empirical elementary scores
$\myS^{\rm Q}_{0.90, \hsp \theta}$ in the bottom row of Figure
\ref{fig:Murphy} confirm this intuition, in line with what we saw in
Figure \ref{fig:wind} for the mean functional.

\section{Discussion}  \label{sec:discussion}

We have studied mixture representations of Choquet type for the
scoring functions that are consistent for quantiles and expectiles,
respectively, including the ubiquitous case of the mean or expectation
functional, and nesting probability forecasts for binary events as a
further special case.  A particularly interesting aspect of these
results is that they allow for an economic interpretation of
consistent scoring functions in terms of betting and investment
problems.  Our interpretation of expectiles as optimal decision
thresholds in investment problems with fixed costs and differential
tax rates appears to be original and may bear on the current debate
about the revision of the Basel protocol for banking regulation.
 
From a general applied perspective, Gneiting (2011, p.~757) had
argued that if point forecasts are to be issued and evaluated, 
\begin{quote} 
\small
``it is essential that either the scoring function be specified ex
ante, or an elicitable target function be named, such as the mean or a
quantile of the predictive distribution, and scoring functions be used
that are consistent for the target functional.''
\end{quote} 
Patton (2015, p.~1) took this argument a step further, by positing
that 
\begin{quote} 
\small
``rather than merely specifying the target functional, which narrows
the set of relevant loss functions only to the class of loss functions
consistent for that functional [\ldots] forecast consumers or survey
designers should specify the single specific loss function that will
be used to evaluate forecasts.''
\end{quote} 
This is a very valid point.  Whenever forecasters are to be
compensated for their efforts in one way or another, the scoring
function ought to be disclosed.  To give an example of this best
practice, the participants of forecast competitions hosted on the
Kaggle platform (\url{www.kaggle.com}) are routinely informed about
the relevant scoring function prior to the start of the
competition. See, e.g., Hong et al.~(2014) for a description of the
Global Energy Forecasting Competition 2012.

However, there remain many situations in which point forecasters
receive directives in the form of a functional, without an
accompanying scoring function being available.  This might be, because
the forecasts are utilized by a myriad of communities, a situation
often faced by national and international weather centers, because
costs and losses are unknown or confidential, because the goal is
general methodological development, as opposed to a specific applied
task, because interest centers on an understanding of forecasters'
behaviors and performance, or simply because of negligence of best
practices.  In such settings, our findings suggest the routine use of
new diagnostic tools in the evaluation and ranking of forecasts, which
we call Murphy diagrams.  Interest sometimes centers on decompositions
of expected or empirical scores into uncertainty, resolution, and
reliability components, as studied by DeGroot and Fienberg (1983),
Br\"ocker (2009), and Bentzien and Friederichs (2014), among others.
Extensions of Murphy diagrams in these directions may be worthwhile.

Our results also bear on estimation problems, in that scoring
functions connect naturally to M-estimation (Huber 1964; Koltchinskii
1997).  An interesting observation is that the loss functions that
have traditionally been employed for estimation in quantile
regression, ordinary least squares regression, and expectile
regression, namely the asymmetric piecewise linear and squared error
scoring functions (\ref{eq:apl}) and (\ref{eq:ase}), correspond to the
choice of the Lebesgue measure in the mixture representations
(\ref{eq:mr4q}) and (\ref{eq:mr4e}), respectively.  This is in
contrast to binary regression, where estimation is typically based on
the logarithmic score, which corresponds to the choice of the infinite
measure with density $h(\theta) = (\theta \hsp (1-\theta))^{-1}$ in
the mixture representation (\ref{eq:mr4p}), rather than the Lebesgue
or uniform measure that yields (half) the Brier score
(\ref{eq:Brier}).  Quite generally, this raises the question of the
optimal choice of the loss or scoring function to be used for
estimation in regression problems.  Focusing on the binary case, Hand
and Vinciotti (2003), Buja et al.~(2005), Lieli and Springborn (2013)
and Elliott et al.~(2015) have considered the use of economically
motivated criteria.  The interpretations developed in the present
paper can help design economically motivated criteria in more general
settings.

Mixture representations of Choquet type can be found for other, more
general classes of consistent scoring functions.  For instance, our
results extend to the class of functionals known as generalized
quantiles or M-quantiles (Breckling and Chambers 1988; Koltchinskii
1997; Bellini et al.~2014; Steinwart et al.~2014), which subsume both
quantiles and expectiles.  Related, but more complex mixture
representations apply in the case of scoring functions that are
consistent for multi-dimensional functionals, as recently studied by
Fissler and Ziegel (2015).

An interesting question is whether there might be mixture
representations in terms of economically interpretable elementary
scores for proper scoring rules.  As noted, a scoring rule
$\hat{\myS}(F,y)$ assigns a loss or penalty when we issue the
predictive CDF $F$ and $y$ realizes, and for a scoring rule to be
proper, the expectation inequalities in (\ref{eq:proper}) need to
hold.  As we have seen, a predictive distribution for a binary
variable can be identified with a probability forecast, so the
representation (\ref{eq:mr4p}) applies and the answer is well known to
be positive in this case.  However, an extension from probability
forecasts of binary to ternary or general discrete variables does not
appear to be feasible, due to results by Johansen (1974) and
Bronshtein (1978) in convex analysis.\footnote{In a nutshell, Savage
(1971) showed that in the case of $k + 1$ categories, the proper
scoring rules for probability forecasts essentially are parameterized
by the convex functions on the unit simplex in $\real^k$.  Johansen
(1974) and Bronshtein (1978) proved that if $k \geq 2$ then the
extremal members of that class lie dense.}  Despite this negative
result, a closer look at a popular score is encouraging.
Specifically, the widely used continuous ranked probability score
(CRPS; Matheson and Winkler 1976), 
\[ 
\hat{\myS}(F,y) =
\int_{-\infty}^\infty \left( F(\theta) - \one(\theta \geq y) \right)^2
\hsp \dd \theta, 
\] 
equals the integral of the Brier score (\ref{eq:Brier}) for the
induced probability forecast, namely $F(\theta)$, of the binary event
$\{ Y \leq \theta \}$ over all thresholds ${\theta} \in \real$.  For
simplicity, let us assume that $F$ has unique quantiles.  We may then
invoke the mixture representation (\ref{eq:mr4p}) along with the
relationships (\ref{eq:psfx}) and (\ref{eq:qep}) to yield\footnote{An
expected or empirical CRPS then corresponds to the volume under the
surface spanned by the Murphy diagrams for all $\alpha$-quantile
predictions, or by the Murphy diagrams for all threshold-determined
binary probability forecasts.}
\[ 
\hat{\myS}(F,y) 
= 2 \int_{-\infty}^{+\infty} \! \int_0^1 \myS^{\rm B}_\alpha( F(\theta),
    \one(\theta \geq y) \, \dd\alpha \, \dd\theta 
= 2 \int_{-\infty}^{+\infty} \! \int_0^1 \mySQ( q_{\alpha,F},y) \,
    \dd\alpha \, \dd\theta.  
\] 
Depending on the order of integration, the mixture representation
recovers the quantile or the threshold decomposition of the CRPS
(Gneiting and Ranjan 2011) after evaluating the first integral.  More
complex weighting schemes depending on $\theta$ and $\alpha$ can be
employed, for a general family of proper scoring rules that can be
economically motivated and justified.  Related ideas have recently
been put forward in the hydrologic and meteorological literatures
(Laio and Tamea 2007; Bradley and Schwartz 2011; Smet et al.~2012).

\section*{Acknowledgement}

This work has been funded by the European Union Seventh Framework
Programme under grant agreement no.~290976.

\section*{References}

\newenvironment{reflist}{\begin{list}{}{\itemsep 0mm \parsep 1mm
\listparindent -7mm \leftmargin 7mm} \item \ }{\end{list}}

\vspace{-6.5mm}
\begin{reflist}

{\rm Bellini, F., Klar, B., M\"uller, A. and Rosazza Gianin,
  E.~R.}~(2014) Generalized quantiles as risk measures. {\em
  Insurance: Mathematics and Economics}, 54, 41--48.

{\rm Bellini, F. and Bignozzi, V.}~(2015) On elicitable risk
measures. {\em Quantitative Finance}, 15, 725--733.

{\rm Bentzien, S. and Friederichs, P.}~(2014) Decomposition and
graphical portrayal of the quantile score. {\em Quarterly Journal of
  the Royal Meteorological Society}, 140, 1924--1934.

{\rm Berrocal, V. J., Raftery, A. E., Gneiting, T. and Steed,
  R. C.}~(2010) Probabilistic weather forecasting for winter road
maintenance. {\em Journal of the American Statistical Association},
105, 522--537.

{\rm Bradley, A.~A. and Schwartz, S.~S.}~(2011) Summary
Verification measures and their interpretation for ensemble
forecasts. {\em Monthly Weather Review}, 139, 3075--3089.

{\rm Breckling, J. and Chambers, R.}~(1988) M-Quantiles. {\em
Biometrika}, 75, 761--771.

{\rm Br\"ocker, J.}~(2009) Reliability, sufficiency, and the
decomposition of proper scores. {\em Quarterly Journal of the Royal
Meteorological Society}, 135, 1512--1519.

{\rm Bronshtein, E.~M.}~(1978) Extremal convex functions. {\em
Siberian Journal of Mathematics}, 19, 6--12.

{\rm Buja, A., Stuetzle, W. and Shen, Y.}~(2005) Loss functions for
binary class probability estimation and classification: Structure and
applications. Working paper, 
\url{http://www-stat.wharton.upenn.edu/~buja/PAPERS/paper-proper-scoring.pdf}.

{\rm DeGroot, M.~H. and Fienberg, S.~E.}~(1983) The comparison and
evaluation of forecasters. {\em Statistician}, 32, 12--22.

{\rm Delbaen, F., Bellini, F., Bignozzi, V. and Ziegel, J.~F.}~(2014)
Risk measures with the CxLS property. Preprint, \url{arXiv:1411.0426}.

{\rm Diebold, F.~X. and Mariano, R.~S.}~(1995) Comparing predictive
accuracy. {\em Journal of Business and Economic Statistics}, 13,
253--263.

{\rm Efron, B.}~(1991) Regression percentiles using asymmetric squared
error Loss. {\em Statistica Sinica}, 1, 93--125.

{\rm Elliott, G., Ghanem, D. and Kr\"uger, F.}~(2015) Forecasting
conditional probabilities of binary outcomes under misspecification.
Working paper, Department of Economics, University of California at
San Diego, \url{https://sites.google.com/site/fk83research/papers}.

{\rm Embrechts, P., Puccetti, G., R\"uschendorf, L., Wang, R. and
Beleraj, A.}~(2014) An academic response to Basel 3.5. {\em
Risks}, 2, 25--48. 

{\rm Engelberg, J., Manski, C.~F. and Williams, J.}~(2009) Comparing
the point predictions and subjective probability distributions of
professional forecasters. {\em Journal of Business and Economic
  Statistics}, 27, 30--41.

{\rm Faust, J. and Wright, J.~H.}~(2013) Forecasting inflation.  In
{\em Handbook of Economic Forecasting} (eds G. Elliott and
A. Timmermann), vol.~2A, pp.~2--56. Amsterdam: Elsevier.

{\rm Fissler, T. and Ziegel, J.~F.}~(2015) Higher order elicitability
and Osband's principle.  Preprint, \url{arXiv:1503.08123}.

{\rm Gneiting, T.}~(2011) Making and evaluating point forecasts.
{\em Journal of the American Statistical Association}, 106, 746--762.

{\rm Gneiting, T. and Katzfuss, M.}~(2014) Probabilistic
forecasting. {\em Annual Review of Statistics and Its Application},
1, 125--151.

{\rm Gneiting, T. and Raftery, A.~E.}~(2007) Strictly proper
scoring Rules, prediction, and estimation. {\em Journal of the
American Statistical Association}, 102, 359--378.

{\rm Gneiting, T. and Ranjan, R.}~(2011) Comparing density forecasts
using threshold and quantile weighted proper scoring rules. {\em
  Journal of Business and Economic Statistics}, 29, 411--422.

{\rm Gneiting, T. and Ranjan, R.}~(2013) Combining predictive
distributions. {\em Electronic Journal of Statistics}, 7,
1747--1782.

{\rm Gneiting, T., Balabdaoui, F. and Raftery, A.~E.}~(2007)
Probabilistic forecasts, calibration and sharpness. {\em Journal of
  the Royal Statistical Society Series B: Statistical Methodology},
69, 243--268.

{\rm Gneiting, T., Larson, K., Westrick, K., Genton, M.~G. and
  Aldrich, E.}~(2006) Calibrated probabilistic forecasting at the
Stateline wind energy center: The re\-gime-switching space-time
method. {\em Journal of the American Statistical Association},
101, 968--979.

{\rm Hand, D. J. and Vinciotti, V.}~(2003) Local versus global Models
for classification problems: Fitting models where it matters.  {\em
  The American Statistician}, 57, 124--131.

{\rm Holzmann, H. and Eulert, M.}~(2014) The role of the information
set for forecasting -- with applications to risk management. {\em
Annals of Applied Statistics}, 8, 595--621.

{\rm Hong, T., Pinson, P. and Fan, S.}~(2014) Global Energy
Forecasting Competition 2012. {\em International Journal of
Forecasting}, 30, 357--363.

{\rm Huber, P.~J.}~(1964) Robust estimation of a location
parameter. {\em Annals of Mathematical Statistics}, 35,
73--101.

{\rm Johansen, S.}~(1974) The extremal convex functions. {\em
  Mathematica Scandinavica}, 34, 61--68.

{\rm Koenker, R.}~(2005) {\em Quantile Regression}. Cambridge: Cambridge
University Press.

{\rm Koenker, R. and Bassett, G.}~(1978) Regression quantiles.  {\em
  Econometrica}, 46, 33--50.

{\rm Koltchinskii, V.~I.}~(1997) M-Estimation, convexity and
quantiles. {\em Annals of Statistics}, 25, 435--477.

{\rm Kr\"amer, W.}~(2005) On the ordering of probability
forecasts. {\em Sankhy$\bar{a}$}, 67, 662--669.

{\rm Laio, F. and S.~Tamea} (2007) Verification tools for
probabilistic forecasts of continuous hydrological variables. {\em
  Hydrology and Earth System Sciences}, 11, 1267--1277.

{\rm Lieli, R.~P. and Springborn, M.}~(2013) Closing the gap between
risk estimation and decision making: Efficient management of
trade-related invasive species risk. {\em Review of Economics and
  Statistics}, 95, 632--645.

{\rm Matheson, J.~E. and Winkler, R.~L.}~(1976) Scoring rules for
continuous probability distributions. {\em Management Science},
22, 1087--1096.

{\rm Merkle, E. and Steyvers, M.}~(2013) Choosing a strictly
proper scoring rule. {\em Decision Analysis}, 10, 292--304.

{\rm Murphy, A. H.}~(1977) The value of climatological, categorical
and probabilistic forecasts in the cost-loss ratio situation. {\em
  Monthly Weather Review}, 105, 803--816.

{\rm Murphy, A.~H. and Winkler, R.~L.}~(1987) A general framework
for forecast verification. {\em Monthly Weather Review}, 115,
1330--1338.

{\rm Mylne, K.~R.}~(2002) Decision making from probability
forecasts based on forecast value. {\em Meteorological
  Applications}, 9, 307--315.

{\rm Newey, W.~K. and Powell, J.~L.}~(1987) Asymmetric least
squares estimation and testing. {\em Econometrica}, 55, 819--847.

{\rm Newey, W.~K. and West, K.~D.}~(1987) A simple, positive
semi-definite, heteroskedasticity and autocorrelation consistent
correlation matrix. {\em Econometrica}, 55, 703--708.

{\rm Patton, A.~J.}~(2011) Volatility forecast comparison using
imperfect volatility proxies. {\em Journal of Econometrics}, 160,
246--256.

{\rm Patton, A.~J.}~(2015) Comparing possibly misspecified
forecasts. Working paper, Department of Economics, Duke University,
\url{http://public.econ.duke.edu/~ap172/Patton_bregman_comparison_27mar15.pdf}.

{\rm Phelps, R.~R.}~(2001) {\em Lectures on Choquet's Theorem}, 2nd edn.
Heidelberg: Springer.

{\rm Richardson, D.~S.}~(2000) Skill and relative economic value of of
the ECMWF ensemble prediction system. {\em Quarterly Journal of the
  Royal Meteorological Society}, 126, 649--667.

{\rm Richardson, D.~S.}~(2012) Economic value and skill. In {\em
  Forecast Verification: A Practitioner's Guide in Atmospheric
  Science} (eds I.~T Jolliffe and D.~B. Stephenson), 2nd edn.,
pp.~167--184. Chichester: Wiley.

{\rm Rudebusch, G.~D. and Williams, J.~C.}~(2009) Forecasting
recessions: The puzzle of the enduring power of the yield curve.  {\em
  Journal of Business and Economic Statistics}, 27, 492--503.

{\rm Savage, L.~J.}~(1971) Elicitation of personal probabilities and
expectations. {\em Journal of the American Statistical Association},
66, 783--810.

{\rm Schervish, M.~J.}~(1989) A general method for comparing
probability assessors. {\em Annals of Statistics}, 17,
1856--1879.

{\rm Schulze Waltrup, L., Sobotka, F., Kneib, T. and Kauermann,
  G.}~(2014) Expectile and quantile regression --- David and Goliath?
{\em Statistical Modelling}, \url{doi:10.1177/1471082X14561155}.

{\rm Shuford Jr, E. H., Albert, A. and Massengill, H. E.}~(1966)
Admissible probability measurement procedures. {\em Psychometrika},
31, 125--145.

{\rm Smet, G., Termonia, P. and Deckmyn, A.}~(2012) Added economic
value of limited-area multi-EPS weather forecasting applications.
{\em Tellus Series~A}, 64, 18901.

{\rm Steinwart, I., C.~Pasin, R.~C.~Williamson and Zhang, S.}~ (2014)
Elicitation and identification of properties. {\em Journal of Machine
  Learning Research}, 35, 1--45.

{\rm Thompson, J. C. and Brier, G. W.} (1955) The economic utility
of weather forecasts. {\em Monthly Weather Review}, 83, 249--253.

{\rm Thomson, W.}~(1979) Eliciting production possibilities from a
well-informed manager. {\em Journal of Economic Theory}, 20,
360--380.

{\rm Tsyplakov, A.}~(2014) Theoretical guidelines for a partially
informed forecast examiner. Working paper,
\url{http://mpra.ub.uni-muenchen.de/55017/}.

{\rm Vardeman, X. and Meeden, G.}~(1983) Calibration, sufficiency, and
dominance relations for Bayesian probability assessors. {\em Journal
  of the American Statistical Association}, 78, 808--816.

{\rm Wilks, D. S.}~(2001) A skill score based on economic value for
probability forecasts. {\em Meteorological Applications}, 8, 209--219.

{\rm Wolfers, J. and Zitzewitz, E.}~(2008) Prediction markets in
theory and practice. In {\em The New Palgrave Dictionary of Economics}
(eds S.~N. Durlauf and L.~E. Blume), 2nd edn. London: Palgrave
McMillan.

{\rm Ziegel, J.~F.}~(2014) Coherence and elicitability. {\em
Mathematical Finance}, \url{DOI: 10.1111/mafi.12080}. 

\end{reflist}

\section*{Appendix A: Proofs} 

The specific structure of the scoring functions in (\ref{eq:qsf}) and
(\ref{eq:esf}) permits us to focus on the case $\alpha = 1/2$ in the
subsequent proofs, with the general case $\alpha \in (0,1)$ then being
immediate.

\subsection*{A1 \: Proof of Theorems 1a and 1b} 

In the case of quantiles, the mixture representation (\ref{eq:mr4q}),
the fact that $\dd H(\theta) = \dd g(\theta)$ for $\theta \in \real$,
and the relationship $H(x) - H(y) = \myS(x,y)/(1-\alpha)$ for $x >
y$, are straightforward consequences of the fact that for every $g \in
\cI$ and $x, y \in \real$,
\[
g(x) - g(y) = 
\int_{-\infty}^{+\infty} \{ \one(\theta < x) - \one(\theta < y) \} \, \dd g(\theta). 
\]
As the increments of $H$ are determined by $\myS$, the mixing measure
is unique.

\medskip
Turning now to the case of expectiles, we associate with any function
$\phi \in \cC$ the Bregman type function of two variables
\begin{equation}  \label{eq:Phi} 
\Phi(x,y) = \phi(y) - \phi(x) - \phi'(x) \hsp (y-x) 
\qquad (x, y \in \real).
\end{equation}
Then the mixture representation (\ref{eq:mr4e}), the fact that $\dd
H(\theta) = \dd \phi'(\theta)$ for $\theta \in \real$, and the
relationship $H(x) - H(y) = \partial_2 \myS(x,y)/(1-\alpha)$ for $x
> y$, are immediate consequences of the fact that for {all} $\phi \in
\cC$ and {$x < y$},
\begin{eqnarray*} 
\Phi(x,y) 
\ceq (y-\theta) \, \phi'(\theta) \bigg|_{\theta=x}^y 
     + \int_x^y \phi'(\theta) \, \dd\theta \\
\ceq \int_x^y (y-\theta) \, \dd\phi'(\theta) 
\, =  \, 2 \int_{-\infty}^{+\infty} \myS^{\rm E}_{1/2,\theta}(x,y) \, \dd\phi'(\theta).
\end{eqnarray*} 
The case $x > y$ is handled analogously, and the case $x = y$ is
trivial.  Finally, as the increments of $H$ are determined by $\myS$,
the mixing measure is unique.

\subsection*{A2 \: Proof of Propositions 1a and 1b} 

In the case of the elementary quantile scoring function
(\ref{eq:esfx}), suppose that $\mySQ = (\myS_1 + \myS_2)/2$, where
$\myS_1$ and $\myS_2$ are of the form (\ref{eq:qsf}) with associated
functions $g_1, g_2 \in \cI_1$.  Then
\[
0 = \left\{
\begin{array}{ll} 
( \hsp g_1(x) - g_1(y)) + ( \hsp g_2(x) - g_2(y)) - 2, & y \leq \theta < x, \\
( \hsp g_1(x) - g_1(y)) + ( \hsp g_2(x) - g_2(y)) + 2, & x \leq \theta < y, \rule{0mm}{4mm} \\ 
( \hsp g_1(x) - g_1(y)) + ( \hsp g_2(x) - g_2(y)),     & \rm{otherwise}. \rule{0mm}{4mm} \\ 
\end{array} \right.
\]
As $g_1, g_2 \in \cI_0$ we have $g_j(x) - g_j(y) \in [0,1]$ if $y \leq
x$, and $g_j(x) - g_j(y) \in [-1,0]$ if $x \leq y$, where $j = 1,
2$. It follows that $g_1(x) - g_1(y) = g_2(x) - g_2(y) = 1$ in the
first case, $g_1(x) - g_1(y) = g_2(x) - g_2(y) = - 1$ in the second
case, and $g_1(x) - g_1(y) = g_2(x) - g_2(y) = 0$ in the third case.
This coincides with the value distribution of $g(x) - g(y)$ when $g(x)
= \one(\theta < x)$, whence indeed $\myS_1 = \myS_2 = \mySQ$.

\medskip
In the case of the elementary expectile scoring function
(\ref{eq:esfx}), suppose that $\mySE = (\myS_1 + \myS_2)/2$, where
$\myS_1$ and $\myS_2$ are of the form (\ref{eq:esf}) with associated
functions $\phi_1, \phi_2 \in \cC_1$.  Let $\Phi_1, \Phi_2$ be defined
as in (\ref{eq:Phi}).  Then
\[
\Phi_1(x,y) + \Phi_2(x,y) - 2 \hsp \myS^{\rm E}_{1/2,\theta}(x,y) = 0.  
\]
Taking left-{hand} derivatives with respect to $y$, we obtain 
\[
\left( \phi_1'(x) - \phi_1'(y) \right) + 
\left( \phi_2'(x) - \phi_2'(y) \right) - 
2 \left( \one(\theta < x) - \one(\theta < y) \right) = 0.
\]
As $\phi_1', \phi_2' \in \cI_1$, we may apply the same argument as in
the quantile case to show that $\phi_1'(x) - \phi_1'(y) = \phi_2'(x) -
\phi_2'(y) = \one(\theta < x) - \one(\theta < y)$, whence $\myS_1 =
\myS_2 = \mySE$.

\subsection*{A3 \: Proof of Propositions 2a and 2b} 

In the case of the elementary quantile scoring function $\mySQ$ in
(\ref{eq:qsfx}) suppose first that $x_2 < x_1 \leq q_{\alpha,F}$.
Since
\[
\mySQ(x_2,y) - \mySQ(x_1,y) = 
\left( \one(y \leq \theta) - \alpha \right) 
\left( \one(\theta < x_2) - \one(\theta < x_1) \right) \! ,
\]
we have 
\[
\myE_F \hsp [\mySQ(x_2,Y)] - \myE_F \hsp [\mySQ(x_1,Y)] 
= \left( F(\theta) - \alpha \right)  
  \left( \one(\theta < x_2) - \one(\theta < x_1) \right) \! .
\]
The second factor on the right-hand side vanishes unless $\theta \in
[x_2, x_1)$, and under this latter condition we have $F(\theta) \leq
  \alpha$ and $\one(\theta < x_2) - \one(\theta < x_1) = - 1$, {whence} the
  desired expectation inequality.  The case $q_{\alpha,F} \leq x_1 <
  x_2$ is handled analogously.

\medskip
In the case of the elementary expectile scoring function $\mySE$ in
(\ref{eq:esfx}) we assume first that $x_2 < x_1 \leq t$, where 
$t$ denotes the $\alpha$-expectile of $F$.  Since
\[
\mySE(x_2,y) - \mySE(x_1,y) = 
\left( (1-\alpha) (\theta- y)_+ -  \alpha (y-\theta)_+ \right) 
\left( \one(\theta < x_2) - \one(\theta < x_1) \right) \! ,
\]
we get 
\begin{eqnarray*} 
\lefteqn{\hspace{-125mm} \myE_F \hsp [\mySE(x_2,Y)] - \myE_F \hsp [\mySE(x_1,Y)]} && \\
= \left( (1-\alpha) \myE_F (\theta- Y)_+ - \alpha \hsp \myE_F (Y-\theta)_+ \right) 
\left( \one(\theta < x_2) - \one(\theta < x_1) \right) \! .
\end{eqnarray*}
As the first term on the right-hand side is strictly increasing in
$\theta$ and has a unique zero at the $\alpha$-expectile of $F$, the
proof can be completed in the same way as above.

\section*{Appendix B: Details for the synthetic example}

\begin{table}[t]

\centering

\footnotesize

\begin{tabular}{lll}
\toprule
Forecast & $\alpha$-Quantile & Mean \rule{0mm}{3mm} \\
\midrule
$F$      & $\myE_\myQ \, \mySQ(q_{\alpha,F},Y)$ & $\myE_\myQ \, \myS^{\rm E}_{1/2, \hsp \theta}(\mu_F,Y)$ \rule{0mm}{3.5mm} \\
\midrule
Perfect 
& $a_{\alpha,\theta} + \alpha \hsp \Phi(\theta - z_\alpha) + \int_{\theta-z_\alpha}^\infty \! A_\theta(x) \hsp \dd x$ 
& $c_\theta - \theta \hsp \Phi(\theta) - \varphi(\theta)$ \rule{0mm}{4mm} \\
Climatological 
& $a_{\alpha,\theta} + \min \! \left( \Phi(\frac{\theta}{\sqrt{2}}), \alpha \right)$ 
& $c_\theta - \theta \hsp \one(\theta \geq 0)$  \rule{0mm}{4mm} \\
Unfocused 
& $a_{\alpha,\theta} 
+ \myE_\tau \!\! \left[ \alpha \hsp \Phi(\theta - z_{\alpha,\tau}) + \int_{\theta-z_{\alpha,\tau}}^\infty A_\theta(x) \hsp \dd x \right]$ 
& $c_\theta - \myE_\tau \!\! \left[ \theta \hsp \Phi \! \left( \theta - \frac{\tau}{2} \right) 
  + \varphi \! \left( \theta - \frac{\tau}{2} \right) \right]$ \rule{0mm}{4mm} \\
Sign-reversed 
& $a_{\alpha,\theta} + \alpha \hsp \Phi(\theta - z_\alpha) + \int_{-\infty}^{z_\alpha - \theta} A_\theta(x) \hsp \dd x$ 
& $c_\theta - \theta \hsp \Phi(\theta) + \varphi(\theta)$ \rule{0mm}{4mm} \\
\bottomrule
\end{tabular}

\caption{\small Expected extremal scores in the prediction space example of Table
  \ref{tab:sim}.  For $\alpha \in (0,1)$ and $\theta \in \real$, we
  let $a_{\alpha,\theta} = - \alpha \hsp \Phi(\theta/\sqrt{2})$,
  $A_\theta(x) = \Phi(\theta-x) \hsp \varphi(x)$, and $c_\theta =
  \theta \hsp \Phi(\theta/\sqrt{2}) + \sqrt{2} \hsp
  \varphi(\theta/\sqrt{2})$, where $\Phi$ and $\varphi$ denote the CDF
  and the probability density function of the standard normal
  distribution, respectively.  
  \label{tab:expscores}}

\end{table}

Here we give details for the synthetic example introduced in Table
\ref{tab:sim} and discussed throughout Section \ref{sec:rankings}.
Table \ref{tab:expscores} shows analytic expressions for the expected
score
\[
\myE_\myQ \, \myS( \hsp \myT(F), Y) 
\] 
where $F$ is either the perfect, the climatological, the unfocused, or
the sign-reversed forecaster, and the functional $\myT(F)$ is either
the $\alpha$-quantile, $q_{\alpha,F}$, or the mean, $\mu_F$, of the
CDF-valued random quantity $F$.  The scoring function $\myS$ is the
elementary quantile scoring function $\mySQ$ in (\ref{eq:qsfx}) or the
elementary scoring function $\myS^{\rm E}_{1/2, \hsp \theta}$ in
(\ref{eq:esfx}).  For example, if $X$ is a quantile forecast for $Y$
at level $\alpha \in (0,1)$ then
\begin{equation}  \label{eq:decom}  
\myE_\myQ \, \mySQ( \hsp X, Y) 
= - \alpha \, \myQ( Y \leq \theta) + \alpha \, \myQ( X \leq \theta) 
  + \myQ( X > \theta, Y \leq \theta). 
\end{equation} 
decomposes into three terms, the first depending on the outcome only,
the second depending on the forecast only, and the third accounting
for the joint distribution.  In view of the relationships
(\ref{eq:psfx}) and (\ref{eq:qep}), the foregoing covers the case of
the extremal scoring function $\mySB$ for event probabilities, too.

\end{document}